\documentclass[a4paper,11pt]{article}

\usepackage{authblk,bbm,fancyhdr}
\usepackage{amsmath,amsthm,amsfonts,amstext,amssymb,amsopn,mathtools}
\usepackage{algorithm,algpseudocode,enumitem,standalone}
\usepackage{tikz}
\usetikzlibrary{shapes, arrows,decorations,automata,backgrounds,petri,positioning,arrows}

\usepackage[a4paper, margin=1in, headheight=13.6pt]{geometry}

\newcommand{\argmax}{\mathrm{argmax}}
\newcommand{\argmin}{\mathrm{argmin}}

\newtheorem{theorem}{Theorem}[section]
\newtheorem{lemma}[theorem]{Lemma}

\theoremstyle{remark}
\newtheorem{remark}{Remark}[section]

\title{Stochastic Dynamic Linear Programming: A Sequential Sampling Algorithm for Multistage Stochastic Linear Programming}
\author[1]{Harsha Gangammanavar\thanks{harsha@smu.edu}}
\author[2]{Suvrajeet Sen\thanks{s.sen@usc.edu}}
\affil[1]{Department of Engineering Management, Information, and Systems, Southern Methodist University, Dallas TX}
\affil[2]{Department of Industrial and Systems Engineering, University of Southern California, Los Angeles, CA}

\date{}
\newcommand{\theoremRef}[1]{{\scshape Theorem} \ref{#1}}
\newcommand{\lemmaRef}[1]{{\scshape Lemma} \ref{#1}}
\newcommand{\assumRef}[1]{(A\ref{#1})}
\newcommand{\algRef}[1]{Algorithm \ref{#1}}

\newcommand{\dynamics}[2]{\mathcal{D}_{#1}(#2)}
\newcommand{\minorants}[2]{\mathcal{J}_{#1}^{#2}}

\newcommand{\child}{\Omega}

\newcommand{\expect}[2]{\mathbb{E}_{#2}[#1]}
\newcommand{\history}[2]{{#1}_{[#2]}}
\newcommand{\future}[2]{{#1}_{(#2)}}
\newcommand{\basis}{\mathbb{B}}
\newcommand{\basisSet}{\mathcal{B}}
\newcommand{\inner}[2]{\langle #1, #2 \rangle}

\begin{document}\thispagestyle{empty}
\maketitle

\begin{abstract}
Multistage stochastic programming deals with operational and planning problems that involve a sequence of decisions over time while responding to realizations that are uncertain. Algorithms designed to address multistage stochastic linear programming (MSLP) problems often rely upon scenario trees to represent the underlying stochastic process. When this process exhibits stagewise independence, sampling-based techniques, particularly the stochastic dual dynamic programming (SDDP) algorithm, have received wide acceptance. However, these sampling-based methods still operate with a deterministic representation of the problem that uses the so-called sample average approximation. In this work, we present a sequential sampling approach for MSLP problems that allows the decision process to assimilate newly sampled data recursively. We refer to this method as the stochastic dynamic linear programming (SDLP) algorithm. Since we use sequential sampling, the algorithm does not necessitate a priori representation of uncertainty, either through a scenario tree or sample average approximation, both of which require a knowledge/estimation of the underlying distribution. In this regard, SDLP is a sequential sampling approach to address MSLP problems. This method constitutes a generalization of the Stochastic Decomposition (SD) for two-stage SLP models.  We employ quadratic regularization for optimization problems in the non-terminal stages. Furthermore, we introduce the notion of basic feasible policies which provide a piecewise-affine solution discovery scheme, that is embedded within the optimization algorithm to identify incumbent solutions used for regularization. Finally, we show that the SDLP algorithm provides a sequence of decisions and corresponding value function estimates along a sequence of state trajectories that asymptotically converge to their optimal counterparts, with probability one.
\end{abstract}

\section{Introduction} \label{sect:intro}
Many practical applications require sequences of decisions to be made under evolving and often uncertain conditions. Multistage stochastic programming (SP) is one of the common approaches used to guide decision making in such stochastic optimization problems.  A variety of fields ranging from traditional production systems \cite{Peters1977}, hydroelectric reservoir scheduling \cite{Morton1996, Pereira1991}, and financial planning models \cite{Carino1998, Kusy1986}, to emerging applications in electricity grids with renewable generation \cite{Powell2012} and revenue management \cite{Topaloglu2009}, among others, have successfully used multistage SP.

Multistage stochastic linear programming (MSLP) models with recourse were used to formulate the early applications of multistage SP. These MSLP models were solved using multistage extensions of the L-shaped method \cite{VanSlyke1969}, such as the Nested Benders Decomposition (NBD) method \cite{Birge1985a}, the scenario decomposition method \cite{Mulvey1995}; and the progressive hedging (PH) algorithm \cite{Rockafellar1991}. A common feature across all these algorithms is the use of approximate deterministic representation of uncertainty through scenario trees (i.e., precedence relations) built using scenario generation techniques (e.g., \cite{Dupacova2000}). When the underlying stochastic process becomes complicated, their deterministic representation may result in large, unwieldy scenario trees. To handle such scenario trees in a computationally viable manner, one may have to resort to scenario reduction methods (e.g., \cite{Dupacova2003}). For models that allow stagewise independent data, \cite{Pereira1991} proposed the stochastic dual dynamic programming (SDDP) algorithm. The multistage extensions of the L-shaped method and SDDP and its variants intend to solve a base model with an uncertainty representation involving a finite sample space and known probability distribution (a scenario tree or a sample average approximation). The model resulting from such a representation is deterministic in nature. In this regard, we refer to these methods as deterministic decomposition-based methods.

In problems where reliable knowledge of uncertainty is not available, an approach that does not rely on exact probabilistic information is desirable. For MSLP models, the first inexact bundle method proposed in \cite{Sen2014} called the multistage stochastic decomposition (MSD) achieves this objective. This algorithm is a dynamic extension of the regularized version of the two-stage stochastic decomposition (2-SD) algorithm \cite{Higle1994}. It accommodates very general stochastic processes, possibly with time correlations, through a nodal formulation which requires only a ``layout" of a scenario tree, and a mechanism that provides transitions between nodes. A standard scenario tree formulation is a special case of such a mechanism. When the stochastic process exhibits interstage independence, a time-staged formulation (as opposed to nodal scenario-tree formulation) is more convenient. With this in mind, we present a sequential sampling-based algorithm that addresses decision making under stagewise independent stochastic processes.

\subsection{Contributions} \label{sect:contributions}
We refer to our sequential sampling-based approach for MSLP with interstage independence as the \emph{stochastic dynamic linear programming} (SDLP) algorithm. In light of the existing deterministic and stochastic decomposition-based methods, the contributions of this work are as follows.

\subsubsection*{An Algorithm for Stagewise Independent MSLP Models}
SDLP harnesses the advantages offered by both the interstage independence of stochastic processes (like SDDP) as well as the sequential sampling design (like 2-SD) to build an algorithm. The algorithm achieves asymptotic convergence while sampling only a small number of scenarios (e.g., one) in any iteration. The algorithm is designed for a state variable formulation of the MSLP models. There are many distinguishing features of SDLP when compared to deterministic decomposition-based methods.  The principal differences are highlighted below.
\begin{itemize}
	\item {\it Static v. Dynamic Instances:} Deterministic decomposition-based methods, including SDDP, can be classified as external sampling methods where the uncertainty representation step precedes the optimization step. In such methods, one begins by first identifying the nodes (observations) and the probability of observing the nodes at each stage, which is then used to set up the MSLP instance. SDDP aims to optimize the resulting MSLP instance. The decisions provided by SDDP are justified using the mathematical theory of sample average approximation. The uncertainty representation (observations and probabilities) is explicitly used in computing the cost-to-go value function approximations. In contrast to that, SDLP accommodates the possibility of observing new scenarios during the course of the algorithm. As a result, the uncertainty representation, and therefore, the MSLP instance dynamically evolves with the introduction of new scenarios.
	\item {\it Implications of Sampling:} SDLP completes the forward and backward recursion computations along a single sample-path that is generated independently of previously observed sample-paths. Although this feature is reminiscent of SDDP variants that incorporate sampling in the forward and backward passes, there are two main differences. (a) Since SDDP operates with a fixed uncertainty representation, the sampled paths selected for forward and backward pass calculations are a subset of sample-paths used in the uncertainty representation. On the other hand, the sample-path used in the forward and backward recursions of SDLP may include observations that have not been encountered before. (b) Since the number of observations increases, the piecewise affine approximations need to be updated to ensure that they continue to provide a lower bound for the dynamically changing sample average approximation.
	\item {\it Asymptotic Behavior:} Unlike SDDP that can recover the cost-to-go value functions in finitely many steps, we show the optimality of SDLP using the primal-dual relationships that are fundamental to mathematical programming. Moreover, SDLP approximations are not finitely convergent. Asymptotic convergence distinguishes the mathematical underpinnings of SDLP and SDDP analyses.
\end{itemize}

The distinguishing features identified above are all consequences of sequential sampling. In the two-stage setting (as in the regularized 2-SD algorithm of \cite{Higle1994}), the recourse function is a deterministic optimization problem, a linear program to be specific. On the other hand, in the multistage setting, the recourse functions in non-terminal stages will dynamically update nested sample average approximations. Therefore, MSD as well as SDLP include provisions to address the stochasticity in value function approximations. Since MSD works with a layout of a scenario tree, it uses a node-specific approximation. With stagewise independent stochastic processes, the future value function approximations are shared by all observations at a stage. Therefore, updates along the current sample-path perturb the future approximations for all observations. This marks a subtle but significant difference in the way the approximations are constructed and updated in SDLP. This also impacts the convergence analysis.

\subsubsection*{A Policy to Identify Incumbent Solutions}
The use of quadratic regularization in two-stage SP algorithms (\cite{Higle1994} and \cite{Ruszczynski1986}) has proven to be very effective for several reasons. The quadratic regularizer helps ensure descent, which is a property that helps prove convergence because it imparts approximate (or estimated) monotonicity.  This property was very useful for convergence proofs, as in \cite{Higle1994}, for two-stage SLP problems.  Another important advantage is that one can limit the size of the stage optimization problem to at most $n_t+3$ ``cuts'', where $n_t$ is the number of decision variables in stage $t$. Motivated by the advantages offered by regularization in sampling-based two-stage algorithms, the proposed algorithm, as well as MSD, employ quadratic regularization. Quadratic regularization can also be interpreted in the context of proximal algorithms at all non-terminal stages using ``incumbent" decisions\footnote{In SP algorithms, especially methods based on 2-SD, an ``incumbent decision" is one for which the predicted objective value appears to be the best (at the current iteration). When predictions change, the incumbent decision must also be updated.} that are maintained for all sample-paths discovered during the algorithm. Maintaining and updating these incumbent solutions becomes cumbersome as the number of sample-paths increases. To address this critical issue we develop the notion of a piecewise-affine policy which is used to identify incumbent solutions for out-of-sample scenarios (new sample-paths) generated sequentially within the algorithm. Such a policy is referred to as a basic feasible policy (BFP). A BFP is based on the optimal bases of the approximate stage problems that are solved during the course of the algorithm. While the BFP designed in this paper is used to identify incumbent solutions for SDLP, the general idea underlying a BFP can also be adopted for other multistage SP algorithms, including SDDP.

This paper also serves as a companion to our earlier work \cite{Gangammanavar2018} by providing the theoretical corroboration of the empirical evidence presented there. In \cite{Gangammanavar2018}, a sequential sampling-based approach was used for controlling distributed storage devices in power systems with significant renewable resources. Computational experiments conducted on large-scale instances showed that such an approach provides solutions which are statistically indistinguishable from solutions obtained using SDDP, while significantly reducing the computational time. These improvements (in comparison to SDDP) can be attributed to two key features of the SDLP algorithm. Firstly, the forward and backward recursion calculations are carried out only along one sample-path. This significantly reduces the total number of optimization problems solved in any iteration. Secondly, the use of regularization allows us to maintain a fixed-sized optimization problem at each stage, as in the case of the master problem in the regularized 2-SD algorithm \cite{Higle1994}. This implies that the computational effort per iteration (necessary to solve stagewise optimization problems) does not increase with iterations. Moreover, it has been recently established that 2-SD  provides a sequence of incumbent solutions that converges to the optimal solution at a sublinear convergence rate \cite{Liu2020asymptotic}.  It is important to emphasize that this result pertains to a solution sequence, rather than the objective function sequence, which was already known for first-order methods such as stochastic approximation (SA).  Because the design and analysis of this paper mirrors that of 2-SD, we suspect that a similar rate of convergence may be possible for SDLP as well. However, a detailed convergence rate analysis is beyond the scope of the current paper.

\subsection*{Organization} The remainder of the paper is organized as follows. In \S\ref{sect:form} we present the MSLP formulation used in this paper. We present a brief overview of the deterministic decomposition-based MSLP methods, particularly SDDP, in \S\ref{sect:mslpAlgorithms}.   A detailed description of the SDLP algorithm is provided in \S\ref{sect:sdlp}. We present the convergence analysis of SDLP in \S\ref{sect:convergenceAnalysis}. Our presentation will have a particular emphasis on the differences in approximations employed in deterministic and stochastic decomposition methods. 

\section{Notation and Formulation} \label{sect:form}
We consider a system where sequential decisions are made at discrete decision epochs denoted by the set $\mathcal{T} := \{0,\ldots,T\}$. Here $T < \infty$, and hence we have a finite horizon sequential decision model with $T+1$ stages. In the interest of brevity (especially because there are many subscripted elements) we denote by $t+$ and $t-$ the succeeding and preceding time periods $(t+1)$ and $(t-1)$, respectively.  We use $[t]$ to denote the history of the stochastic process $\{v_t\}_{t=0}^T$ until (and including) stage $t$, i.e., $v_{[t]} = v_0, v_1,\ldots,v_t$. Likewise, we use $v_{(t+)}$ to denote the process starting from stage $t+1$ until the end of horizon (stage $T$), i.e., $v_{(t+)} = v_{t+1},\ldots,v_T$. We use $\inner{\cdot}{\cdot}$ to denote the inner product of vectors (e.g., $\inner{v}{w} = v^\top w$) and the product of a matrix transpose and a vector, i.e., $\inner{M}{v} = M^\top v$.

Commonly in SP, MSLP models are formulated without state variables, focusing only on decisions in each stage.  However in many applications, especially those involving dynamic systems, it is common to use the state variable description of system evolution.  Because we expect SDLP to be able to provide decision support for such systems, it is advisable to use a state variable formulation.  This approach is also common in the dynamic programming community. In this regard, we use a state variable $s_t := (x_t, \omega_t) \in \mathcal{S}_t$ to describe the system at stage $t$. This state variable is comprised of two components: $x_t \in \mathcal{X}_t$ is the endogenous state of the system and $\omega_t \in \Omega_t$ captures the exogenous information revealed in interval $(t-1,t]$. A stochastic process over which the decision-maker cannot exert any control drives the exogenous state evolution. For example, the exogenous state variable may represent a weather phenomenon like wind speed, or a market phenomenon like the price of gasoline. The evolution of the endogenous state, on the other hand, can be controlled by an algorithm through decisions $u_t$ and is captured by stochastic linear dynamics:
\begin{align} \label{eq:stDyn}
	x_{t+} = \dynamics{t+}{x_t,\omega_{t+}, u_t} = a_{t+} + A_{t+} x_t + B_{t+}u_t.
\end{align}
Here, $(a_{t+}, A_{t+}, B_{t+})$ are components of the exogenous information vector $\omega_{t+}$ corresponding to the next time period.

To characterize the exogenous process $\{\tilde{\omega}_t\}_{t=1}^T$, we use $(\Omega, \mathcal{F}, \mathbb{P})$ to denote the filtered probability space. Here, $\Omega = \Omega_1 \times \ldots \times \Omega_T$ denotes the set of outcomes and $\omega_t$ denotes an observation of the random variable $\tilde{\omega}_t$. The $\sigma$-algebras $\mathcal{F}_t \subseteq \mathcal{F}$ represent the data available to the decision-maker at time $t$, which satisfy $\mathcal{F}_t \subseteq \mathcal{F}_{t^\prime}$ for $t < t^\prime$. The exogenous data $\omega_t$ includes components of $(a_t, A_t, B_t)$ that appear in \eqref{eq:stDyn} and parameters $(b_t, C_t)$ in the right-hand side of the constraints at stage $t$.

With these notations, the state-variable representation of the time-staged MSLP model can be written in the nested form as follows:
\begin{align} \label{eq:mslp}
    \min~& \inner{c_0}{x_0} + \inner{d_0}{u_0} + \mathbb{E}_{\tilde{\omega}_{(1)}}\Bigg[\inner{c_1}{x_1} + \inner{d_1}{u_1} + \mathbb{E}_{\tilde{\omega}_{(2)}} \bigg [\ldots + \\& \hspace{7cm}\mathbb{E}_{\tilde{\omega}_T}  [\inner{c_T}{x_T} + \inner{d_T}{u_T}] \bigg] \Bigg] \notag \\
    \text{s.t.}~& u_t \in \mathcal{U}_t(s_t) := \{u_t|D_tu_t \leq b_t - C_t x_t,~ u_t \geq 0\} \qquad \forall t \in \mathcal{T} \notag \\ 
    & x_{t+} = \dynamics{t+}{x_t,\omega_{t+}, u_t} = a_{t+} + A_{t+} x_t + B_{t+}u_t \qquad \forall t \in \mathcal{T} \setminus \{T\}. \notag
\end{align}
The above problem is stated for a given initial endogenous state $x_0$. Here, $u_t$ for $t \in \mathcal{T}$ are decision vectors and $\mathcal{U}_t(s_t)$ are closed convex sets that define the feasible set of decisions. In our finite horizon framework, we assume that the terminal cost $h_{T+}(s_{T+})$ is known for all $s_{T+}$ (or negligible enough to be set to $0$). The expectation is taken with respect to the exogenous stochastic process $\tilde{\omega}_{(t)}$ over the remainder of the horizon. In a time period $t$, the state $s_t$ explicitly depends on the initial state $x_0$, past decisions $u_{[t]}$, and past exogenous states $\omega_{[t]}$. Since $s_t$ affects the feasible set $\mathcal{U}_t$, decision $u_t$ is the function of decision process until time $t$. The multistage program can alternatively be stated in the following recursive form for all $t \in \mathcal{T}$:
\begin{align}\label{eq:mslpt}
	h_t(s_t) = \inner{c_t}{x_t} + \min~& \inner{d_t}{u_t} + \expect{h_{t+}(\tilde{s}_{t+})}{} \\
	s.t.~& u_t \in \mathcal{U}_t(s_t) := \{u_t|D_tu_t \leq b_t - C_t x_t, u_t \geq 0\}, \notag
\end{align}
where $\tilde{x}_{t+} = \dynamics{t+}{x_t,\tilde{\omega}_{t+},u_t}$.  Since the initial state $x_0$ is assumed to be given, the stage-$0$ (henceforth known as the root-stage) problem has deterministic input.

In general, the MSLP problems are PSPACE-hard \cite{Dyer2006, Hanasusanto2016comment} and require exponential effort in horizon $T$ for provably tight approximations with high probability. To keep our presentation consistent with our algorithmic goals, we make the following assumptions:
\begin{enumerate}
	\renewcommand{\labelenumi}{(A\theenumi)}
	\item The set of root-stage decisions $\mathcal{U}_0$ is compact. \label{assum:compact}
	\item The complete-recourse assumption is satisfied at all non-root stages, that is, the feasible set $\mathcal{U}_t(s_t)$ is non-empty for all state trajectories $s_t$ with $x_t$ satisfying \eqref{eq:stDyn} for all $t \in \mathcal{T} \setminus \{0\}$. \label{assum:completeResource}
	\item The constraint matrices $D_t$ are fixed and have full row rank.\label{assum:fixed}
	\item Zero provides the lower bound on all cost-to-go value functions. \label{assum:zeroLB}
	\item The stochastic process for exogenous information is stagewise independent and its support is finite. \label{assum:indep}
\end{enumerate}
These assumptions provide a special structure and are fairly standard in the SP literature (\cite{Philpott2008, Sen2014}). The fixed recourse assumption \assumRef{assum:fixed} implies that the recourse matrix $D_t$ does not depend on exogenous information. As for assumption \assumRef{assum:zeroLB}, note that most loss functions used in engineering applications and statistical learning obey this property.  For situations in which this assumption is not satisfied, one can perform a pre-processing step as follows: first estimate a lower bound on the optimal objective function value for each stage, and then, add the absolute value of the most negative stagewise lower bound to all stages.  Introducing such a constant into the objective function does not alter the optimal decisions while rendering the validity of \assumRef{assum:zeroLB}. The finite support assumption \assumRef{assum:indep} on exogenous information ensures that $\mathcal{F}_t$ is finite. We note that the algorithms presented here can be extended, after some refinement, to settings where some of the above assumptions can be relaxed. For instance, certain extensions to Markovian stochastic processes can be envisioned. However, a detailed treatment of these extensions is beyond the scope of this paper.

\section{MSLP Algorithms} \label{sect:mslpAlgorithms}
The fundamental difficulty of solving SP problems is associated with the \emph {nested multidimensional} integral for computing the expectation in \eqref{eq:mslpt}. The most direct approach involves incorporating simulation to estimate the expected recourse function as:
\begin{align} \label{eq:saa}
	\expect{h_{t+}(\tilde{s}_{t+})} \approx \widehat{H}_{t+}^N(s_{t+}) := \frac{1}{N}\sum_{n=1}^N h_{t+}(s_{t+}^n)
\end{align}
where, $s_{t+}^n$ has components $x_{t+}^n = a_{t+}^n + A_{t+}^n x_t + B_{t+}^n u_t$ and $\omega_{t+}^n$. Doing so results in the so-called sample average approximation (SAA) problem. In this case, we can view the support of $\Omega_{t+}$ as consisting of a simulated sample $\Omega_{t+}^N := \{\omega_{t+}^1, \omega_{t+}^2,\ldots,\omega_{t+}^N\}$, where each observation vector $\omega_{t+}^n$ has the same probability $p(\omega_{t+}^n) = (1/N)~\forall n = 1,\ldots,N$. Since the recourse function in \eqref{eq:mslpt} involves the expectation operator, it is worth noting that the estimate in \eqref{eq:saa} is an unbiased estimator and under certain conditions (e.g., when the sample is independent and identically distributed) a consistent estimator of the expected recourse function. However, the optimal value of a SAA problem provides a downward biased estimator of the true optimal value \cite{Shapiro2011}.

The SAA problem can be reformulated as a single large linear program (the deterministic equivalent form \cite{Birge2011}), and off-the-shelf optimization software can be used to solve the problem. However, as the sample size increases (as mandated by SAA theory to achieve high-quality solutions \cite{Shapiro2014}), or the number of stages increases, such an approach becomes computationally burdensome. Deterministic decomposition-based cutting plane methods, also known as outer-linearization methods, provide a means to partially overcome the aforementioned burden.

The deterministic decomposition-based (DD) methods can be traced to Kelley \cite{Kelley1960} for smooth convex optimization problems, Benders decomposition for ideas of decomposition/ partitioning in mixed-integer programs (MIPs) \cite{Benders1962}, and Van Slyke and Wets for 2-SLPs \cite{VanSlyke1969}. While the exact motivation for these methods arose in different contexts, we now see them as being very closely related to the outer-linearization perspective. These ideas have become the main-stay for both 2-SLPs and stochastic MIPs. 

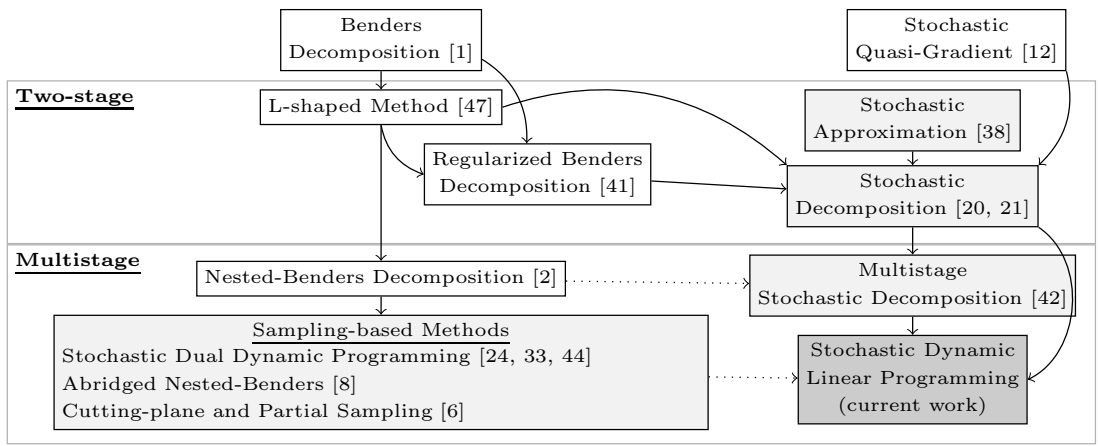
\begin{figure}[t]
\centering
\fontsize{6}{7.2}\selectfont
\begin{tikzpicture}[->,alg/.style={rectangle,draw,align=center}]
	\node[draw=black!30, minimum height=1.8cm, minimum width=0.75\textwidth] (bg2) at (0,0) {};
	\node[alg,above right=0.1cm and 3cm of bg2.north west] (bd) {Benders \\ Decomposition \cite{Benders1962}};	
	\node[alg,below= 0.2cm of bd] (l) {L-shaped Method \cite{VanSlyke1969}}; 
	\node[alg,below right=0.2cm and 1.8cm of l.south west] (rbd) {Regularized Benders \\ Decomposition \cite{Ruszczynski1986}};

	\node[alg,right=4cm of bd] (sqg) {Stochastic \\ Quasi-Gradient \cite{Ermoliev1983}};		
	\node[alg,fill=gray!10, below left =0.2cm and 0.5cm of sqg.south east]  (sa) {Stochastic \\ Approximation \cite{Robbins1951}};
	\node[alg,fill=gray!10, below=0.15cm of sa] (sd) {Stochastic \\ Decomposition \cite{Higle1991,Higle1994}};
	\node[anchor=north west, align=center] at (bg2.north west) {\tiny \bf \underline{Two-stage}};
	
	\draw (bd) to (l); \draw[bend left] (bd) to (rbd); \draw[bend right] (l.south) to (rbd.west); 
	\draw [bend left] (l.east) to (sd.north west);	\draw (rbd) to (sd);  
	 \draw (sa) to (sd);  \draw[bend left] (sqg.south east) to (sd.north east); 
	
	\node[draw=black!30,minimum height=2.2cm, minimum width=0.75\textwidth, below=0cm of bg2] (bg3) {};
	\node[anchor=north west] at (bg3.north west)  {\tiny \bf \underline{Multistage}};
	
	\node[alg,below=1.5cm of l] (nbd) {Nested-Benders Decomposition \cite{Birge1985a}};	
	\node[alg,fill=gray!10, below=0.2cm of nbd] (sbd) {\underline{Sampling-based Methods}\\
	\begin{minipage}{0.44\textwidth}	 
		Stochastic Dual Dynamic Programming \cite{Infanger1996, Pereira1991,Shapiro2011}\\
		Abridged Nested-Benders \cite{Donohue2006}\\
		Cutting-plane and Partial Sampling \cite{Chen1999}		
	\end{minipage}};	
	\node[alg,fill=gray!10, below=0.3cm of sd] (msd) {Multistage \\ Stochastic Decomposition \cite{Sen2014}};
	\node[alg,fill=gray!40, below=0.2cm of msd] (sdlp) {Stochastic Dynamic \\ Linear Programming \\ (current work)};
	
	\draw (l) to (nbd); \draw (nbd) to (sbd);
	\draw (sd) to (msd); 
	\draw[bend left=60] (sd.south east) to (sdlp.east);
	\draw (msd) to (sdlp);
	\draw[dotted] (nbd) to (msd); \draw[dotted] (sbd) to (sdlp);
	
\end{tikzpicture}
\caption{Multistage Stochastic Linear Programming Algorithms} \label{fig:multistageAlgs}
\end{figure}

DD-based algorithms originally developed for 2-SLP have been extended to successive stages of dynamic linear programs. One of the early successes was reported in \cite{Birge1985a}, where the classical two-stage Benders decomposition algorithm was extended to multiple stages. This procedure has subsequently come to be known as the NBD algorithm. The starting point of this algorithm is the scenario tree representation of underlying uncertainty where all possible outcomes and their interdependence are represented as nodes on a tree. Naturally, this implies that the NBD algorithm can be classified under the multistage DD-based methods. Relationships between various algorithmic approaches are summarized in Figure \ref{fig:multistageAlgs}.

\subsection{Stochastic Dual Dynamic Programming}
It is well known that the number of nodes in the scenario tree grows exponentially with the number of stages, and therefore, the need to visit all the nodes in the scenario tree significantly increases the computational requirements of the NBD algorithm. Pereira and Pinto \cite{Pereira1991} provided a sampling-based approach to address this issue in the stochastic dual dynamic programming (SDDP) algorithm.
 
Like MSLP algorithms mentioned earlier, SDDP creates an outer approximation of the stage value function using subgradient information. SDDP performs its iteration in a forward pass and a backward pass, a feature common to most multistage SP algorithms. However, it avoids the intractability of scenario trees by assuming that the stochastic process is stagewise independent. While the algorithm traverses forward using sampling, the approximations are created on the backward pass similar to deterministic Benders type cuts. The interstage independence assumption allows these cuts to be shared across different states within a stage. Cut sharing under special stagewise dependency is presented in \cite{Infanger1996}, the algorithmic enhancements proposed in \cite{Linowsky2005}, and the inclusion of risk measures \cite{Guigues2012sampling, Philpott2012} have extended the capabilities of the original algorithm \cite{Pereira1991} and have contributed to the success of SDDP. The abridged nested decomposition algorithm in \cite{Donohue2006} and the cutting plane and partial sampling algorithm proposed in \cite{Chen1999} are other sampling-based methods which are similar in flavor to SDDP. 

\sloppy
The main steps of SDDP are presented in \algRef{alg:sddp}. As in the case of NBD, each iteration of SDDP begins by solving an optimization problem for the root-stage. Then a finite number of Monte Carlo simulations are carried out to identify forward sample-paths $\{\future{\omega}{0}\}_{n=1}^N$ for the iteration. Along each one of these sample-paths, the forward pass involves identifying candidate solutions $u_t^{kn}$ by solving an optimization problem of the form:
\begin{align} \label{eq:sddpt}
	\min \{f_t^{k-1}(s_t, u_t)~|~u_t \in \mathcal{U}_t(s_t^{kn})\}
\end{align}
and propagating the state according to the dynamics in \eqref{eq:stDyn} as $x_{t+}^{kn} = \mathcal{D}_{t+}(x_t^{kn},\omega_{t+}^{kn}, u_t^{kn})$. These two steps are undertaken in an alternating manner for all stages until the end of the horizon. In the above stage optimization problem, $f_t^{k-1}(s_t, u_t)$ denotes the current approximation of the cost-to-go value function in \eqref{eq:mslpt}. At the end of the forward pass, we have a set of candidate solutions at each non-terminal stage $\{u_t^{kn}\}_{\forall t}$; one for each simulated sample-path of the forward pass.

\begin{algorithm}[!ht]
\caption{Stochastic Dual Dynamic Programming} \label{alg:sddp}
\begin{algorithmic}[1]
\State \textbf{Initialization}: Iteration count $k \leftarrow 0$.
\State {\bf Forward pass:} Decision simulation along simulated sample-paths.\label{alg:sddpLoop}
\State Solve the root-stage optimization problem \eqref{eq:sddpt} to identify $u_0^k$.
\State Sample a set of $N$ paths $\{\future{\omega}{0}^{kn}\}_{n=1}^N$.
	\For {$t = 1,\ldots,T-1$}
		\For {$n = 1,\ldots,N$}
			\State Setup the candidate states $x_t^{kn} = \dynamics{t}{x_{t-}^{kn}, \omega_t^{kn}, u_{t-}^{kn}}$.
			\State Solve the stage optimization problem in \eqref{eq:sddpt} with $s_t^{kn}$ as input, and obtain the optimal primal candidate solution $u_t^{kn}$.
		\EndFor
	\EndFor
\State {\bf Backward pass}: Update cost-to-go value function approximations.
	\For {$t = T-1,\ldots,0$} \label{alg:backwardpassBegin}
		\For {$n = 1,\ldots,N$}
			\For{$\omega_{t+} \in \Omega_{t+}$} \label{alg:2dd_omegaSelection}
			\State Setup $s_{t+} = (x_{t+},\omega_{t+})$, where $x_{t+} = \mathcal{D}_{t+}(x_t^{kn}, \omega_{t+}, u_t^{kn})$.
			\State Solve subproblem with $s_{t+}$ as input:
			\begin{align}\label{eq:2dd_subproblem}
				\min~\{f_{t+}^k(s_{t+},u_{t+})~|~ u_{t+} \in \mathcal{U}_{t+}(s_{t+})\},
			\end{align}
			\hspace*{45pt} and obtain optimal dual solution $\pi_{t+}(\omega_{t+})$.
			\State Compute lower bounding affine function $\ell_{t+}(s_{t+}) := \alpha_{t+}^{kn}(\omega_{t+}) + \inner{\beta_{t+}^{kn}(\omega_{t+})}{x_{t+}}$, where
			\begin{align} \label{eq:2dd_coeff}
				\alpha_{t+}^{kn}(\omega_{t+}) = \inner{b_{t+}}{\pi_{t+}(\omega_{t+})}; \quad \beta_{t+}^{kn}(\omega_{t+}) = c_{t+} - \inner{C_{t+}}{\pi_{t+}(\omega_{t+})}.
			\end{align}
			\State Update the set of coefficients as: $$\minorants{t+}{k}(\omega_{t+}) = \minorants{t+}{k-1}(\omega_{t+}) \cup \{(\alpha_{t+}^{kn}(\omega_{t+}), \beta_{t+}^{kn}(\omega_{t+})\}.$$
			\EndFor		
		\EndFor
		\State Obtain the updated stage cost-to-go value function approximation using		
		\begin{align}\label{eq:sddpApprox}
			h_{t+}^k(s_{t+}) = \max_{j \in \minorants{t+}{k}(\omega_{t+})} \{\alpha_{t+}^j + \inner{\beta_{t+}^j}{x_{t+}}\}.
		\end{align}
		to obtain $f_{t}^k(s_{t}, u_{t}) = \inner{c_t}{x_t} + \inner{d_t}{u_t} + \sum_{\omega_{t+} \in \Omega_{t+}} p(\omega_{t+}) h_{t+}^k(s_{t+})$.
	\EndFor \label{alg:backwardpassEnd}
\State Increment iteration count: $k \leftarrow k + 1$, and go to Line \ref{alg:sddpLoop}.
\end{algorithmic}
\end{algorithm}

In the work of Pereira and Pinto \cite{Pereira1991} the backward pass proceeds as in the case of NBD (see Steps \ref{alg:backwardpassBegin}--\ref{alg:backwardpassEnd} in \algRef{alg:sddp}). At a non-terminal stage $t$ and for each element of the candidate solution set $\{u_t^{kn}\}$, backward pass states are computed using the linear dynamics in \eqref{eq:stDyn} for all possible outcomes in $\Omega_{t+}$. With each of these backward pass states as input, an optimization problem is solved in stage $t+$ and the optimal dual solution is used to compute a lower bounding affine function. Since this procedure requires subproblems to be solved for all the nodes along all the sample-paths simulated in the forward pass, this approach is ideal for narrow trees (few possible realizations per stage). However, the computational issues resurface when the number of outcomes per stage increases. Donohue and Birge proposed the abridged NBD algorithm to address this issue in \cite{Donohue2006} where the forward pass proceeds only along a subset of candidate states (termed as ``branching'' states) while solving all the nodes only along the trajectory of branching states in the backward pass. Subsequently, it was proposed in \cite{Linowsky2005} and \cite{Philpott2008} that sampling procedures can be adopted in the backward pass as well. We make the following observations regarding the original SDDP procedure and its variants:
\begin{enumerate}
	\item Each collection of affine function $\minorants{t}{kn}(\omega)$ is associated with a unique candidate solution $u_{t-}^{kn}$ at stage $(t-1)$. The cost-to-go value function approximations in \eqref{eq:sddpApprox} includes a piecewise linear approximation in which the pointwise maximum is defined over the collections of affine functions generated across all the sample-paths, i.e., $\mathcal{J}_{t+}^k(\omega) = \cup_{n=1}^N \mathcal{J}_{t+}^{kn}(\omega)$. In addition, if the uncertainty is confined to the state dynamics, then the cuts can be shared across the outcomes $\omega \in \Omega_{t+}$. This	``sharing'' of cuts is possible due to the stagewise independence of exogenous information and was first proposed in \cite{Infanger1996}.
	\item A SAA of the problem in \eqref{eq:mslpt} can be constructed by replacing the true distribution of $\tilde{\omega}_t$ by the empirical distribution based on a random sample $\{\omega_t^1, \omega_t^2, \ldots, \omega_t^N\}$ for all $t \in \mathcal{T}\setminus\{0\}$. These random samples are generated independently to ensure that the stagewise independent assumption is respected. A SAA based SDDP algorithm was analyzed in \cite{Shapiro2011}.
	\item The forward pass sampling must ensure that each of the $|\Omega_1| \times |\Omega_2| \times \ldots \times |\Omega_T|$ possible sample-paths are visited infinitely many times w.p.1. If sampling is employed in the cut generation procedure (as in \cite{Chen1999,Philpott2008}), it must be performed independently of the forward pass sampling and must ensure that each element of $\Omega_t$ is sampled infinitely many times w.p.1. at all stages. 
\end{enumerate}

In contrast to the SDDP algorithm, where a fixed sample is used at each stage, our SDLP algorithm will generate approximations that are based on sample average functions constructed using a sample whose size increases with iterations. The sequential nature of introducing new observations into the sample requires additional care within the algorithm design, particularly in the backward pass when approximations are generated (Steps \ref{alg:backwardpassBegin}--\ref{alg:backwardpassEnd}). We present these details in the next section.

\section{Stochastic Dynamic Linear Programming} \label{sect:sdlp}
An iteration of SDLP involves two principal steps: forward and backward recursion. The use of forward and backward recursions is common to almost all the multistage SP algorithms (except those based on progressive hedging \cite{Rockafellar1991}). The forward-backward recursion approach to solving dynamic optimization problems can be traced back to the differential dynamic programming (DDP) algorithm \cite{Jacobson1970differential}. The SDLP algorithm is closely related to the DDP algorithm, in the sense that we create locally accurate approximations of the subdifferential, whereas DDP works with quadratic approximations of smooth deterministic dynamic control problems. The algorithmic constructs of SDLP are designed to accommodate the inherent non-smoothness of MSLP models, and of course,  stochasticity. We present details of these in iteration $k$ of the algorithm. Note that we make the same assumptions as the SDDP algorithm.

\subsection{Forward Recursion} \label{sect:SDLPforwardPass}
The forward recursion begins by solving the following quadratic regularized optimization problem:
\begin{align} \label{eq:regObjFnApprox}
	\min_{u_0 \in \mathcal{U}_0}~\bigg\{f_0^{k-1}(s_0,u_0) + \frac{\sigma}{2}\|u_0 - \hat{u}_0^{k-1}\|^2 \bigg \}.
\end{align}
Here, the proximal parameter $\sigma \geq 1$ is assumed to be given. We denote the optimal solution of the above problem as $u_0^k$ and refer to it as the candidate solution. The incumbent solution $\hat{u}_0^k$ used in the proximal term is similar to that used in the regularized L-shaped \cite{Ruszczynski1986} and 2-SD \cite{Higle1994} algorithms. This is followed by simulating a sample-path $\future{\omega}{0}^k$ that is generated independently of previously observed sample-paths. The remainder of the forward recursion computations is carried out only along this simulated sample-path in two passes - a prediction pass and an optimization pass. 

\subsubsection*{Prediction Pass} At all non-terminal stages we use a regularized stage optimization problem which is centered around the incumbent solution. The goal of the prediction pass is to make sure that the incumbent solutions, and the corresponding incumbent states, satisfy the underlying model dynamics in \eqref{eq:stDyn} along the current sample-path $\future{\omega}{0}^k$. Given the initial state $x_0$, the prediction pass starts by using the root-stage incumbent solution $\hat{u}_0^k$ and computing the incumbent state for stage-$1$ as: $\hat{x}_1^k = \dynamics{1}{x_0, \omega_{1}^k, \hat{u}_0^k}$. At the subsequent stage, we use the BFP to identify the incumbent solutions as $\hat{u}_t^k(\hat{s}_t^k) = \mathcal{M}_t(\hat{s}_t^k)$. Here, $\mathcal{M}_t: \mathcal{S}_t \rightarrow \mathcal{U}_t$ is a vector valued mapping that takes the state vector $s_t$ as an input and maps it on to a solution in $\mathcal{U}_t(s_t)$. We postpone the details of specifying BFP to \S\ref{sect:incumbentSelection} and continue with the algorithm description here. We proceed by computing the incumbent state using \eqref{eq:stDyn} and identifying the incumbent solution using the BFP for the remainder of the horizon. At the end of the prediction pass, we have an  incumbent state $\{\hat{x}_t^k\}$ and solution $\{\hat{u}_t^k\}$ trajectories\footnote{We will use the more explicit notation $\hat{u}_t^k(\hat{s}_t^k)$ that shows the dependence of the incumbent solution on the input incumbent state only when it does not add undue notational burden. In most cases, we will simply use $\hat{u}_t^k$ for incumbent solution.} that satisfy state dynamics in \eqref{eq:stDyn} over the entire horizon. 

\subsubsection*{Optimization Pass} After completing the prediction pass, the optimization pass is carried out to simulate candidate solutions along the current sample-path $\future{\omega}{0}^k$ for all $t \in \mathcal{T}\setminus\{T\}$:
\begin{align}\label{eq:forwardt}
	u_t^k \in \argmin \{ f_t^{k-1}(s_t^k, u_t) + \frac{\sigma}{2}\|u_t - \hat{u}_t^{k}(\hat{s}_t^k) \|^2~|~u_t \in \mathcal{U}_t(s_t^k)\}. 
\end{align}
Here $f_t^{k-1}(s_t,u_t)$ is the current approximation of the cost-to-go value function and the proximal term $\sigma > 0$ is assumed to be given. Structurally, $f_t^{k-1}$ is a piecewise affine and convex function and is similar to the approximations used in the SDDP algorithm. However, each individual piece is a minorant\footnote{Since the approximations generated in sequential sampling-based methods are based on statistical estimates which are updated iteratively, we use the term ``minorant'' to refer the lower bounding affine functions. This usage follows its introduction in \cite{Sen2014} and is intended to distinguish them from the more traditional ``cuts'' in DD-based methods.} generated using certain sample average functions. The candidate decision for a particular stage is used to set up the subsequent endogenous state $x_{t+}^k = \mathcal{D}_{t+}(x_t^k, \omega_{t+}^k, u_t^k)$ and thus the input state $s_{t+}^k$. We refer to the decision problem in \eqref{eq:forwardt} as \textit{Timestaged Decision Simulation} at stage $t$ (TDS$_t$). This completes the optimization pass, and hence the forward recursion, for the current iteration. At the end of forward recursion, we have the incumbent trajectory $\{\hat{x}_t^k\}$ and the candidate trajectory $\{x_t^k\}$ which will be used for updates during the backward recursion.

\subsection{Backward Recursion} \label{sect:SDLPbackwardPass}
The primary goal in the backward recursion procedure is to update the cost-to-go value function approximations $f_t^{k-1}$ at all non-terminal stages. As the name suggests these calculations are carried out backward in time, starting from the terminal stage to the root-stage, along the same sample-path that was observed during the forward recursion. These calculations are carried out for both the candidate as well as the incumbent trajectories. 

In both the DD and SD-based approaches, the value function is approximated by the pointwise maximum of affine functions. However, the principal difference between these approaches lies in how the expected value function is approximated. In DD-based methods, it is the true expected value function which requires the knowledge of the probability distribution or a SAA with a fixed sample (as in \eqref{eq:saa}). On the other hand, the SD-based methods create successive approximations $\{f_t^k\}$ (for $t <T$) that provide a lower bound on a sample average approximation using only $k$ observations in iteration $k$, and therefore, satisfies:
\begin{align} \label{eq:2sdLB}
	f_t^k(s_t, u_t) - \inner{c_t}{x_t} - \inner{d_t}{u_t} \leq \widehat{H}_{t+}^k(s_{t+}) := \sum_{\omega_{t+} \in \Omega_{t+}^k} p^k(\omega_{t+}) h_{t+}(x_{t+},\omega_{t+}),
\end{align}
where $x_{t+}$ is the endogenous state obtained from $\eqref{eq:stDyn}$ with $(x_t,\omega_{t+},u_t)$ as input, for all $\omega_{t+} \in \Omega_{t+}^k$ and $u_t \in \mathcal{U}_t(s_t)$. The quantity $p^k(\omega_{t+})$ in \eqref{eq:2sdLB} measures the relative frequency of an observation which is defined as the number of times $\omega_{t+}$ is observed ($\kappa^k(\omega_{t+})$) over the number of iterations ($k$). This quantity approximates the unconditional probability of exogenous information at stage $t+$, and is updated as follows. Given the current sample-path $\future{\omega}{0}^k$, a collection of observations at a non-root stage $t$ $(t > 0)$ is updated to include the latest observation $\omega_t^k$ as: $\Omega_t^{k} = \Omega_t^{k-1} \cup \omega_t^k$. The observation count is also updated as: $\kappa^k(\omega_t) = \kappa^{k-1}(\omega_t) + \mathbbm{1}_{\omega_t = \omega_t^k}$, for all $\omega_t \in \Omega_t^k$. Using these counts, the observation frequency for $\omega_t \in \Omega_t^k$ is given by: $p^k(\omega_t) = \frac{\kappa^k(\omega_t)}{k}$. Notice the superscript $k$ (iteration count) that is used in our notation of the SAA function $\hat{H}_{t+}^k$, the collection of observations $\Omega_{t+}^k$, and the observation frequency $p^k$. This is intended to convey the sequential nature of SDLP. 

\begin{figure}
	\centering
		\fontsize{6}{7.2}\selectfont		

\begin{tikzpicture}
	
	\node[circle, draw, thick, blue, fill = blue!20, minimum size = 0.4cm] at (0,0) (root) {};
	
	\def\samplPath {(2,0),(2,1),(4,0),(4,-1),(6,0),(8,2),(8,1),(8,-1),(8,-2),(10,0),(10,1),(12,-1)};
	\foreach \observ [count=\i] in \samplPath {
		\node[circle, draw, thick, black!20, minimum size = 0.4cm] at \observ {1};
	};
	\def\samplPath {(2,-1),(4,2),(10,2)};
	\foreach \observ [count=\i] in \samplPath {
		\node[circle, draw, thick, black!20, minimum size = 0.4cm] at \observ {2};
	};

	\node[circle, draw, thick, green, minimum size = 0.4cm, text=red] at (2,1) {2};
	\node[circle, draw, thick, green, minimum size = 0.4cm, text=red] at (4,1) {1};
	\node[circle, draw, thick, green, minimum size = 0.4cm, text=red] at (6,2) {4};
	\node[circle, draw, thick, green, minimum size = 0.4cm, text=red] at (8,0) {1};
	\node[circle, draw, thick, green, minimum size = 0.4cm, text=red] at (10,-1) {1};
	\node[circle, draw, thick, green, minimum size = 0.4cm, text=red] at (12,1) {4};			
	\draw[green] (0.2,0) -- (1.8,1); \draw[green] (2.2,1) -- (3.8,1); \draw[green] (4.2,1) -- (5.8,2);
	\draw[green] (6.2,2) -- (7.8,0); \draw[green] (8.2,0) -- (9.8,-1); \draw[green] (10.2,-1) -- (11.8,1);	

	\node at (0,-2.5) {Root node};
	\foreach \n in {1,2,3,4,5,6} {						
		\node at (2*\n,-2.5) {$\Omega_{\pgfmathprintnumber \n}^5$};		
	}	

	\draw[black, rounded corners] (-1,-3) rectangle (12.5,2.5); 
			
\end{tikzpicture}
	\caption{Uncertainty representation after $5$ iteration. The green path denotes the sample-path observed in iteration $5$. The number on the nodes represent the number of times the observation was encountered, i.e., $\kappa^k(\omega)$. New nodes are added to the representation as and when they are encountered. For example, the second node in $\Omega_2^5$ was encountered for the first time in iteration 5.}\label{fig:sequentialSampling}%
\end{figure}
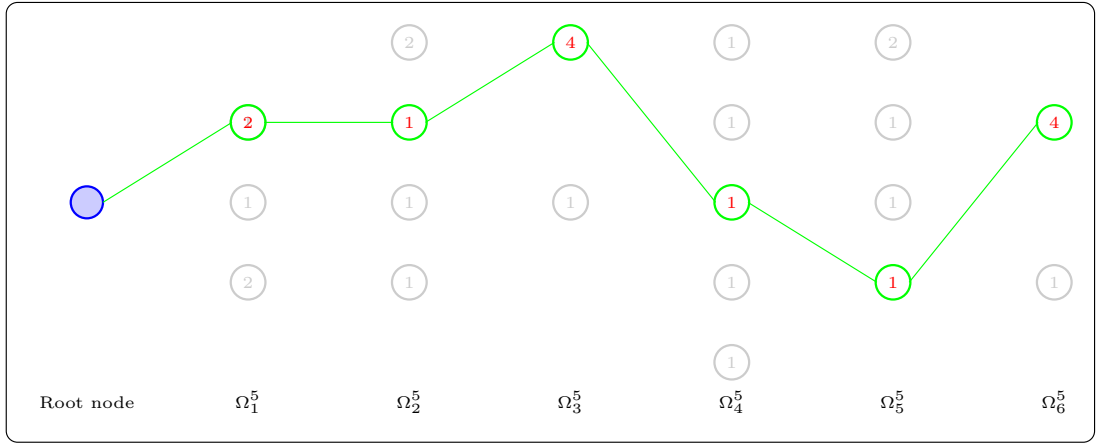

\subsubsection{Terminal Stage Approximation} At the terminal stage, recall that $\expect{h_{T+}(s_{T+})}{} = 0$, and the value function $h_T$ is the value of a deterministic linear program for a given state input $s_T$. The sample average $\widehat{H}_T^k(s_{T}) = \sum_{\omega_T \in \child_T^k} p^k(\omega_T) h_T(s_{T})$ provides an unbiased estimate of $\expect{h_T(\tilde{s}_T)}{}$. Hence, the value function at the penultimate stage ($t = T-1$) can be approximated using a procedure similar to the one employed in the 2-SD algorithm. 

In this procedure, a subproblem corresponding to the current observation $\omega_T^k$ is setup and solved. This subproblem uses $s_T^k = (x_T^k,\omega_T^k)$ as input, where $x_T^k = \mathcal{D}_T(x_{T-1}, \omega_T^k, u_{T-1}^k)$. Let the optimal dual solution obtained be denoted as $\pi_T^k(\omega_T^k)$ which is added to the collection of previously discovered dual vertices: $\Pi_T^k = \Pi_T^{k-1} \cup \pi_T^k(\omega_T^k)$. For other observations in $\Omega_T^k$, i.e., $\omega_T \in \Omega_T^k$ and $\omega_T \neq \omega_T^k$, we identify the best dual vertex in $\Pi_T^k$ using the ``argmax'' operation as in the case of 2-SD algorithm \cite{Higle1991}. This operation is as follows:
\begin{align} \label{eq:argmax}
	\pi_T^k(\omega_T) \in \argmax \{\inner{\pi_T}{(b_T - C_Tx_T)~|~\pi_T \in \Pi_T^k}\}.
\end{align}
Using the dual vertices $\{\pi_T^k(\omega_T)\}_{\omega_T \in \Omega_T^k}$, we  compute the lower bounding affine function $\ell_T^k(s_T) := \alpha_T^k(\omega_T) + \inner{\beta_T^k(\omega_T)}{x_T}$, where
\begin{align} \label{eq:2sd_coeff}
	\alpha_T^k(\omega_T) = \inner{b_T}{\pi_T^k(\omega_T)}; \quad \beta_T^k(\omega_T) = c_T - \inner{C_T}{\pi_T^k(\omega_T)}.
\end{align}
The above calculations are also carried out for the the incumbent state $\hat{s}_t^k$, resulting in the affine function $\hat{\ell}_T^k(s_T) = \hat{\alpha}_T^k(\omega_T) + \inner{\hat{\beta}_T^k(\omega_T)}{x_T}$. The set of affine functions thus obtained  ($\minorants{T}{k} = \minorants{T}{k-1} \cup \{\ell_T^k(s_T), \hat{\ell}_T^k(s_T)\}$) provides the piecewise affine lower bounding function to the value function $h_T(s_T)$ that is given by:
\begin{align} \label{eq:minoT}
	h_T^k(s_T) = \max_{j \in \minorants{T}{k}(\omega_T)} \{\ell_T^j(s_T) = \alpha_T^j(\omega_T) + \inner{\beta_T^j(\omega_T)}{x_T}\}.
\end{align}
The above function provides an outer linearization of the terminal value function.

\subsubsection{Non-terminal Stage Approximation}
When updating the approximations at a non-terminal stage $t$, we have access to the minorants at stage $t+$ (recall that the value functions are being updated recursively backwards from the terminal stage). Using these we can define:
\begin{align} \label{eq:samplMean}
	H_t^k(s_t) := \inner{c_t}{x_t} + \min_{u_t \in \mathcal{U}_t(s_t)}~& \inner{d_t}{u_t} + \sum_{\omega_{t+} \in \child_{t+}^k} p^k(\omega_{t+}) h_{t+}^k(s_{t+}),
\end{align}
where $s_{t+} = (\dynamics{t+}{x_t, \omega_{t+}, u_t}, \omega_{t+})$ for all $\omega_{t+} \in \child_{t+}^k$. The expression in  (\ref{eq:samplMean}) represents a sample average computed over the current observations $\Omega_{t+}^k$ at stage $t+$ at an arbitrary input state $s_t$. Since we use lower bounding approximations $h_{t+}^k$ in building this sample average, this sampled estimate is biased. The stage approximation is updated using a lower bound to the above sample average function, and hence, is biased as well.

In order to compute this lower bound, notice that we can obtain the subgradient, i.e., $\beta_{t+}^k(\omega_{t+}) \in \partial h_{t+}^k(\dynamics{t+}{x_t^k, \omega_{t+}, u_t^k},\omega_{t+})$ using the collection of affine functions $\minorants{t+}{k}(\omega)$ for all observations $\omega_{t+} \in \Omega_{t+}$ (see \S\ref{sect:subgradientPolicy} for details). Let $\alpha_{t+}^k(\omega_{t+})$ be the corresponding intercept term. Using these, a valid lower bound to the sample average function in \eqref{eq:samplMean} can be written as:
\begin{align} \label{eq:samplMeanPrimal}
	H_t^k(s_t) \geq  \inner{c_t}{x_t} + \min_{u_t \in \mathcal{U}_t(s_t)}~\inner{d_t}{u_t} + \sum_{\omega_{t+} \in \child_{t+}^k} p^k(\omega_{t+}) \bigg[\alpha_{t+}^k(\omega_{t+}) + \inner{\beta_{t+}^k(\omega_{t+})}{x_{t+}}\bigg].
\end{align}
Substituting the state dynamics equation in \eqref{eq:stDyn}, and dualizing the linear program on the right-hand side of the above inequality, we obtain:
\begin{align}\label{eq:samplMeanDual}
	H_t^k(s_t) \geq \inner{c_t}{x_t} + &\bar{\alpha}_{t+}^k + \inner{\bar{\beta}_{t+}^k}{x_t} + \\ &\max~\{ \inner{\pi_t}{(b_t - C_tx_t)}~|~\inner{D_t}{\pi_t} \leq \bar{\rho}_{t+}^k,~ \pi_t \leq 0\}, \notag
\end{align}
where, 
\begin{align*}	
	&\bar{\beta}_{t+}^k = \sum_{\omega_{t+} \in \child_{t+}^k} p^k(\omega_{t+}) \inner{\beta_{t+}^k(\omega_{t+})}{A_{t+}},~~ 
	\bar{\rho}_{t+}^k = d_t + \sum_{\omega_{t+} \in \child_{t+}^k} p^k(\omega_{t+}) \inner{\beta_{t+}^k(\omega_{t+})}{B_{t+}}, \\
	&\text{and } \bar{\alpha}_{t+}^k = \sum_{\omega_{t+} \in \child_{t+}^k} p^k(\omega_{t+}) [\alpha_{t+}^k(\omega_{t+}) + \inner{\beta_{t+}^k(\omega_{t+})}{a_{t+}}].
\end{align*}
We refer to the linear program on the right-hand side of inequality in \eqref{eq:samplMeanDual} as the {\it stagewise-dual approximation} at stage $t$ and denote it as  (SDA$_t^k$). Let $\pi_t^k(\omega_t^k)$ denote the optimal dual solution obtained by solving (SDA$_t^k$) with $s_t^k$ as input. Using this we obtain a lower bounding affine function $\ell_t^k(s_t) = \alpha_t^k(\omega_t^k) + \inner{\beta_t^k(\omega_t^k)}{x_t}$ with the following coefficients:
\begin{align}\label{eq:coefft}
	\alpha_t^k(\omega_t^k) =~\inner{\pi_t^k(\omega_t^k)}{b_t} + \bar{\alpha}_{t+}^k~; \qquad \beta_t^k(\omega_t^k) =~ c_t-\inner{C_t}{\pi_t^k(\omega_t^k)} + \bar{\beta}_{t+}^k.
\end{align}
Similar calculations using $\hat{\pi}_t^k(\omega_t^k)$, an optimal solution to the (SDA$_t^k$) with $\hat{s}_t^k$ as input, yields an incumbent affine function $\hat{\ell}_t^k(s_t)$. As before these functions are included in a collection of affine functions to obtain the updated set $\minorants{t}{k}(\omega_t^k)$. 

\begin{algorithm}[!t]
\caption{Stochastic Dynamic Linear Programming} \label{alg:sdlp}
\begin{algorithmic}[1]
\State \textbf{Initialization}: 
\State Choose a proximal parameter $ \sigma \in [\sigma^{min}, \sigma^{max}]$ with $1 \leq \sigma^{min} < \sigma^{max}$.
\State Set observations $\Omega_t^0 = \emptyset$; a trivial affine functions $\ell_t^0 = 0$ in the set $\minorants{t}{0}$ for all $t \in \mathcal{T}$; iteration counter $k \leftarrow 1$.
\State {\bf Forward recursion}: Decision simulation along simulated sample-path \label{alg:sdlpLoop}
\State Solve the root-stage optimization problem of the form \eqref{eq:regObjFnApprox} to identify $u_0^k$.
\State Simulate a sample-path $\future{\omega}{0}^k$.
\State {\it Prediction pass}: 
\For {$t = 1,\ldots,T-1$}
	\State Setup the incumbent state $\hat{x}_t^k = \dynamics{t}{\hat{x}_{t-}^k, \omega_t^k, \hat{u}_{t-}^k(\hat{s}_{t-}^k)}$. 
	\State Identify an incumbent solution $\hat{u}_t^k(\hat{s}_t^k) = \mathcal{M}_t^k(\hat{s}_t^k)$.
\EndFor
\State {\it Optimization pass}:
\For {$t = 1,\ldots,T$}
	\State Setup the candidate state $x_t^k = \dynamics{t}{x_{t-}^k, \omega_t^k, u_{t-}^k}$.
	\State Solve the stage optimization problem \eqref{eq:forwardt} using $s_t^k$ as input, and obtain the \hspace*{0.4cm} candidate primal solution $u_t^k$.
\EndFor
\State {\bf Backward recursion}: Update value function approximations.
\For {$t = T,\ldots,1$}
	\State Setup the stagewise-dual approximation \eqref{eq:samplMeanDual}.
	\State Solve the dual approximation using the candidate and the incumbent states,  \hspace*{0.4cm} and compute the coefficients for affine functions using \eqref{eq:coefft}.
	\State Obtain the updated value function approximation as in \eqref{eq:objUpdtt}.
\EndFor
\State Increment the iteration count $k \leftarrow k + 1$, and go to Line-\ref{alg:sdlpLoop}.
\end{algorithmic}
\end{algorithm}

While it is true that the latest affine functions satisfy $H_t^k(s_t) \geq \ell_t^k(s_t)$, the same does not hold for affine functions generated at earlier iterations. Hence, it is possible that there exists a $j \in \minorants{t}{k}(\omega_t)$ such that the affine function $\ell_t^j(s_t)$ may not lower bound the current sample average $H_t^k(s_t)$. In keeping with the updates of 2-SD \cite{Higle1991}, the old minorants need to be updated as the sample average estimate changes during the course of the algorithm. Under assumption \assumRef{assum:zeroLB}, this is achieved by scaling down the previously generated affine functions. In the two-stage case, 2-SD minorants are updated by multiplying the coefficients by $(k-1)/k$. In the multistage case, the minorants are updated\footnote{The exponent $(T-t)$ results from the fact that minorants in the future $T-t$ stages are also updated in a similar manner. \theoremRef{thm:outerLinearization} provides the formal argument.} as follows
\begin{align}\label{eq:minot}
	h_t^k (s_t)= \max~ \bigg \{ \bigg \{ \bigg (\frac{k-1}{k} \bigg)^{T-t}~\ell_t^j(s_t) \bigg \}_{j \in \minorants{t}{k-1}(\omega_t)},~\ell_t^k(s_t), ~\hat{\ell}_t^k(s_t) \bigg\}. 
\end{align}
Notice that both the candidate and incumbent affine functions generated in previous iterations are treated similarly while scaling down. 

We use these updated minorants to obtain the stage objective function as follows:
\begin{align}\label{eq:objUpdtt}
	f_t^k(s_t, u_t) =~& \inner{c_t}{x_t} + \inner{d_t}{u_t} + \sum_{\omega_{t+} \in \child_{t+}^k} p^k(\omega_{t+}) h_{t+}^k(s_{t+}),
\end{align}
where $s_{t+} = (\dynamics{t+}{x_t, \omega_{t+}, u_t}, \omega_{t+})$, for all $\omega_{t+} \in \child_{t+}^k$. Similar updates are carried out at all the non-terminal stages by progressing backwards to the root-stage along the same sample-path that was used in the forward recursion. The backward recursion for iteration $k$ is said to be complete once the root-stage objective function is updated. The sequentially ordered steps of SDLP algorithm are presented in \algRef{alg:sdlp}. 

\subsubsection{Comparison of DD and SD-based approximations}
The complete recourse assumption ensures that the dual feasible set is non-empty and the optimal dual solution $\pi_T$ is an extreme point of $\{\pi_T ~|~ \inner{D_t}{\pi_t} \leq d_T, \pi_T \leq 0 \}$. There are finitely many of these extreme points, and hence, coefficients for the terminal stage computed using \eqref{eq:2dd_coeff} for the DD-based algorithms or \eqref{eq:2sd_coeff} for the SD-based methods take finitely many values. 

In DD-based multistage algorithms the coefficients belong to a finite set at stage $t+$, and therefore, there exists an iteration $k^\prime$ such that the set of coefficients $\mathcal{J}_{t+}^k(\omega) = \mathcal{J}_{t+}^{k^\prime}(\omega)~\forall \omega \in \Omega_{t+}$ for $k > k^\prime$. Consequently, the dual feasible region of the problem solved in the backward pass has the following form: \begin{align} \label{eq:sddpDualSet}
	\Pi_t^{k,DD} = \left\{ (\pi_t, \theta_{t+}) ~\Bigg \vert~
	\begin{array}{l}
		\inner{D_t}{\pi_t} \leq d_t + \sum_{\omega \in \Omega_{t+}} \sum_{j \in \mathcal{J}_{t+}^k(\omega)} \theta_{t+}^j(\omega) \beta_{t+}^j(\omega), \\
		\sum_{j \in \mathcal{J}_{t+}^k(\omega)} \theta_{t+}^j(\omega) = p(\omega) \quad \forall \omega \in \Omega_{t+},~ \pi_t \leq 0
	\end{array} \right \}.
\end{align}
Notice that this dual feasible region does not change for iterations $k > k^\prime$. Since there are finite number of extreme points to $\Pi_t^{k,DD}$, the coefficients computed using these extreme point solutions result in at most a finite number of distinct values at stage $t$.

In SDLP, notice the update of the old affine functions in \eqref{eq:minot} at stage $t+$ can be viewed as a convex combination of the coefficient vector $(\alpha_{t+}^j, \beta_{t+}^j)$ and a zero vector. Due to these updates, the dual feasible region depends on updated coefficients (particularly $\beta_{t+}^k(\omega)$) as well as frequencies $p^k(\omega)$: 
\begin{align}\label{eq:sdlpDualSet}
\Pi_t^{k,SD} = \{ \pi_t ~\vert~
	\begin{array}{l}
		\inner{D_t}{\pi_t} \leq d_t + \sum_{\omega_{t+} \in \child_{t+}^k} p^k(\omega_{t+}) \inner{\beta_{t+}^k(\omega_{t+})}{B_{t+}}, 
		\pi_t \leq 0
	\end{array} \}.
\end{align}
This implies that dual solutions used to compute the coefficients no longer belong to a finite set. However, following assumption \assumRef{assum:completeResource} the dual feasible set in (SDA$_t^k$) is bounded. Therefore, the coefficients computed in \eqref{eq:coefft} for a non-terminal stage are only guaranteed to be in a compact set. Proceeding backwards, we can conclude that this is the case for coefficients at all non-terminal stages. These observations are summarized in the following lemma.
\begin{lemma} \label{lemma:coeff}
Suppose the algorithm runs for infinitely many iterations. Under assumption Assumption \assumRef{assum:compact} and \assumRef{assum:completeResource}. For all $k \geq 1$,
\begin{enumerate}[label=(\roman*)]
	\item The coefficients of cuts generated within DD-based methods in \eqref{eq:2dd_coeff}, and coefficients of minorants generated for the terminal stage within SD-based methods in \eqref{eq:2sd_coeff} belong to finite sets. \label{lemma:coeff_a}
	\item The coefficients of minorants generated within SD-based methods for the non-terminal stages in \eqref{eq:coefft} belong to compact sets for all $k \geq 1$. \label{lemma:coeff_b}
\end{enumerate}
\end{lemma}
As a consequence of \ref{lemma:coeff_a} in above lemma and \assumRef{assum:indep}, a finite number of cuts are generated during the course of DD-based algorithms for MSLP models. This is possible because these algorithms utilize the knowledge of transition probabilities in computing cut coefficients. Additionally, these cuts provide lower bound to the true value function and are not required to be updated over the course of the algorithm. It must be noted that, the finite number of cuts pertains only to DD-based methods applied to MSLP problems. In the case of multistage stochastic non-linear convex programs (e.g., \cite{Girardeau2015convergence, GuiguesRegularized2020}), the number of cuts is not guaranteed to be finite. In such cases, the coefficients in the non-terminal stages of DD-based methods also belong to compact sets, albeit for a different reason than in the SD-based methods for MSLP models.

The subgradients computed in the SD-based methods are stochastic in nature. Therefore, only affine functions generated in the current iteration satisfy the lower bounding property for the current sample average approximation, but not necessarily for the true value function. The previous affine functions have to be updated using the scheme described in \eqref{eq:minot}. This scheme ensures that the minorant $h_t^k$, obtained after computing the current affine function and updating all previous affine functions, provides a lower bound to the sample average function $H_t^k$ at all non-terminal stages. The outer linearization property of the minorants is formalized in the following theorem.

\begin{theorem}\label{thm:outerLinearization} 
	Suppose assumption \assumRef{assum:compact}-\assumRef{assum:indep} hold. 
	\begin{enumerate}[label=(\roman*)]
		\item The minorant computed in \eqref{eq:minoT} for terminal stage satisfies:
		\begin{subequations}
		\begin{align} \label{eq:outerLinearization_T}
			 h_T(s_T) \geq h_T^k(s_T) \geq h_T^{k-1}(s_T) \geq \ldots \geq h_T^j(s_T),
		\end{align}
		 for all $1 \leq j \leq k$, $s_T \in \mathcal{S}_T$.
		\item At non-terminal stages, the minorant computed in \eqref{eq:minot} satisfies for $s_t \in \mathcal{S}_t$:
	\begin{align} \label{eq:outerLinearization_t}
		H_t^k(s_t) \geq h_t^k(s_t) \geq \bigg(\frac{k-1}{k}\bigg)^{T-t} h_t^{k-1}(s_t).
	\end{align}
	\end{subequations}
	\end{enumerate}
\end{theorem}
\begin{proof}
The first part of the theorem follows directly from the linear programming duality and the construction of the affine functions $\ell_T^k$ in \eqref{eq:argmax} and \eqref{eq:2sd_coeff}. For proof of the second part, we use $m = \omega_{t+}^k$ which is the observation encountered at stage $t+$ in iteration-$k$ and $n$ to index the set $\Omega_{t+}$. Following this notation, we denote $x_{t+} = \dynamics{t+}{x_t, \omega_{t+}, u_t}$ as $x_{nt+}$ and $s_{t+} = (x_{t+}, \omega_{t+})$ as $s_{nt+}$. Consider the stage sample average problem in \eqref{eq:samplMean}:
\begin{align}\label{eq:outer_1}
	H_t^k(s_t) - & \inner{c_t}{x_t} =~ \min_{u_t \in \mathcal{U}_t(s_t)} \inner{d_t}{u_t} + \sum_{n \in \child_{t+}^k} p^k(n) h_{t+}^k(s_{nt+}).
\end{align}
Recall that the affine function $\ell_t^k$ is computed using the dual solution of the problem on the right-hand side of the above equation. Using \eqref{eq:samplMeanDual} and linear programming duality, we obtain
\begin{align}
	H_t^k(s_t) \geq \ell_t^k(s_t).
\end{align}
We distribute the summation in \eqref{eq:outer_1} over observations encountered in the first $i < k$ iterations (i.e., $\Omega_{t+}^i$) and those encountered after iteration $i$.
\begin{align*}
	H_t^k(s_t) - & \inner{c_t}{x_t} \\
			   =~ &\min_{u_t \in \mathcal{U}_t(s_t)} \inner{d_t}{u_t} + 
			   \sum_{n \in \child_{t+}^j} p^k(n) h_{t+}^k(s_{nt+}) + \sum_{n \in \child_{t+}^k \setminus \child_{t+}^j} p^k(n) h_{t+}^k(s_{nt+}). \notag
\end{align*}
 Since $h_{t+}^k \geq 0$, we have 
\begin{align*}			   
	H_t^k(s_t) - \inner{c_t}{x_t} \geq~& \min_{u_t \in \mathcal{U}_t(s_t)} \inner{d_t}{u_t} + \sum_{n \in \child_{t+}^i} p^k(n) h_{t+}^k(s_{nt+}) \\
	=~ & \min_{u_t \in \mathcal{U}_t(s_t)} \inner{d_t}{u_t} + \sum_{n \in \child_{t+}^j} \frac{\kappa^{i}(n) + \kappa^{[i,k]}(n)}{k} \cdot h_{t+}^k(s_{nt+}). 
\end{align*}
For observations in $\Omega_{t+}^i$, we distribute the computation of their relative frequency by setting $\kappa^k(n) = \kappa^i(n) + \kappa^{[i, k]}(n)$, where $\kappa^{[i,k]}(n)$ is the number of times observation $n$ was encountered after iteration $i$. Once again invoking $h_{t+}^k \geq 0$ we obtain:
\begin{align*}
	H_t^k(s_t) - \inner{c_t}{x_t} \geq~ & \min_{u_t \in \mathcal{U}_t(s_t)} \inner{d_t}{u_t} + \sum_{n \in \child_{t+}^i} \frac{i}{k} \times \frac{\kappa^{i}(n)}{i} \cdot h_{t+}^k(s_{nt+}).
\end{align*}
Recall that the minorants at stage $t+$ are updated in \eqref{eq:minot} by adding new affine function into the collection while multiplying the previously generated affine function by a factor of $(\frac{i}{k})^{T-t-1} < 1$. By replacing the current minorant $h_{t+}^k$ by the scaled version of the one available in iteration $j$, we have:
\begin{align}			   
	H_t^k(s_t) \geq ~ & \inner{c_t}{x_t} + \min_{u_t \in \mathcal{U}_t(s_t)} \inner{d_t}{u_t} + \frac{i}{k} \sum_{n \in \child_{t+}^i} p^i(n) \bigg[\bigg(\frac{i}{k}\bigg)^{T-t-1} h_{t+}^i(s_{nt+})\bigg] \notag \\
	\geq~ & \bigg(\frac{i}{k}\bigg)^{T-t} \bigg [ \inner{c_t}{x_t} + \min_{u_t \in \mathcal{U}_t(s_t)} \inner{d_t}{u_t} + \sum_{n \in \child_{t+}^j} p^i(n) h_{t+}^i(s_{nt+}) \bigg]. \notag
\end{align}
The second inequality follows from assumption \assumRef{assum:zeroLB}. Notice that the scaling factor used when $t+ = T$ reduces to one. In this case, the future cost corresponds to the terminal stage, and the affine functions satisfy $\ell_T^j(s_T) \leq h_T(s_T)$ for all $j \in \minorants{T}{k}(\omega_T)$. Therefore, $h_T^k(s_T) \leq h_T(s_T)$. At other stages, an affine function generated in iteration $i < k$, viz. $\ell_t^j(s_t)$ with $j \in \minorants{t}{i}$ provides a lower bound to the sample average in the same iteration $H_t^i(s_t)$. This leads us to conclude that
\begin{align} \label{eq:lowerBound_scaling}
	H_t^k(s_t) \geq \bigg(\frac{i}{k}\bigg)^{T-t} H_t^i(s_t) \geq \bigg(\frac{i}{k}\bigg)^{T-t} \ell_t^j(s_t).
\end{align}
Applying the same arguments for all $i < k$, and using the definition of minorant in \eqref{eq:minot} we obtain $H_t^k(s_t) \geq h_t^k(s_t)$.

Since,
\begin{align*}
    H_t^k(s_t) \geq \bigg(\frac{i}{k}\bigg)^{T-t} \ell_t^i(s_t) =~& \bigg(\frac{k-1}{k}\bigg)^{T-t} \times \bigg(\frac{k-2}{k-1}\bigg)^{T-t} \times \ldots \times \bigg(\frac{i}{i+1}\bigg)^{T-t} \ell_t^i(s_t) \\
    =~& \bigg(\frac{k-1}{k}\bigg)^{T-t} \bigg(\frac{i}{k-1}\bigg)^{T-t} \ell_t^i(s_t) \\ =~& \bigg(\frac{k-1}{k}\bigg)^{T-t} h_t^{k-1}(s_t).
\end{align*}
This completes the proof.
\end{proof}

As noted in the above proof, the scaling factor $(\frac{i}{k})^{T-t}$ used in \eqref{eq:lowerBound_scaling} is applied to affine functions in $\minorants{t}{i}$ that were generated in iteration $i < k$. Since these affine functions are updated in every iteration, computational efficiency can be attained by using recursive updates. In iteration $k$, the affine functions in $\minorants{t}{k-1}$ are updated by multiplying them by the factor $(\frac{k-1}{k})^{T-t}$ and storing the updated minorants in $\minorants{t}{k}$. We refer the reader to \cite{Higle1996} and \cite{Gangammanavar2020sd} for details regarding efficient implementation of these updates. In the next result we capture the asymptotic behavior of the sequence of minorants $\{h_t^k\}$.
\begin{theorem}\label{thm:minorantAsymptotics}
Under assumption \assumRef{assum:completeResource}, \assumRef{assum:fixed} and \assumRef{assum:indep}, the sequence of functions $\{h_t^k\}_k$ is  equicontinuous and uniformly convergent at all non-root stages.
\end{theorem}
\proof Recall that the coefficients of the minorants belong to a compact set at all the non-root stages (\lemmaRef{lemma:coeff}). Therefore, $\{h_t^k\}$ is a sequence of bounded continuous functions with a uniform Lipschitz constant, say $M$. Further, the sequence $\{h_t^k\}$ converges pointwise on $s_t \in \mathcal{S}_t$. Let $s_t^{k_1}$ and $s_t^{k_2}$  be input states such that $\|s_t^{k_1} - s_t^{k_2}\| < \epsilon/M$, for a positive constant $\epsilon$. From Lipschitz continuity, we have 
\begin{align*}
    |h_t^k(s_t^{n_1}) - h_t^k(s_t^{n_2})| \leq M\|s_t^{n_1} - s_t^{n_2}\| < \epsilon.
\end{align*}
for any $k \geq 1$. Hence, the sequence $\{h_t^k\}$ is equicontinuous. Equicontinuity and pointwise convergence together imply uniform convergence \cite{Rudin1976}.
\endproof

In contrast to the above results, the approximations created in the DD-based methods (see \eqref{eq:sddpApprox}) provide outer linearization for a fixed cost function $H_t^N(\cdot)$. Since, the probability distribution is explicitly used (as constants) in computing the DD-based cuts, the approximations improve monotonically over iterations. That is, $H_t^N(s_t) \geq h_t^k(s_t) \geq h_t^{k-1}(s_t)$ for all $s_t$, without any need for updates. We close this section with the following two remarks. The first contrasts the incorporation of sampling during backward recursion of SDDP with the role of sampling adopted in SDLP. The second identifies the online sampling feature of SDLP that has many advantages in practical settings.

\remark{Sampling during backward recursion has also been explored in SDDP(e.g., \cite{Chen1999} \cite{DeMatos2015improving}, and \cite{Philpott2008}). However, there are important factors that distinguish value function updates undertaken during the backward recursion of SDLP when compared to SDDP calculations. In SDLP, the latest sample-path along which the backward recursion calculations are carried out is included independently of previously encountered sample-paths. As a result, the set of sample-paths grow in size (by at most one) when compared to the set of sample-paths used in the previous iteration. If the latest sample-path was not encountered before, it was not included in calculations carried out in the backward recursion of any previous iterations. This is unlike SDDP where the set of sample-paths is fixed and backward pass calculations are carried out over all scenarios in every iteration. Even when sampling is employed in the backward pass of SDDP, calculations are carried out along all sample-paths by either solving a subproblem or using the ``argmax'' procedure in \eqref{eq:argmax}. This type of cut formation was first suggested in \cite{Higle1991}. Even if the latest path was encountered in earlier iterations, the repeated observation results in an update in the empirical frequency associated with nodes along the latest sample-path. As a consequence, the weights (that are synonymous with estimated probability) used in calculating the SDLP cut coefficients \eqref{eq:coefft} differ from one iteration to the next. In SDDP, on the other hand, actual observation probabilities are used to calculate the value function approximation (see \eqref{eq:sddpApprox}) even when sampling is used on the backward pass.} \label{remark:samplingSDLP}

\remark{Since the SDLP algorithm works with data discovered through sequential sampling, it does not rely on any a priori knowledge of exogenous probability distribution. This feature makes this algorithm suitable to work with external simulators or statistical models that can better capture the nature of exogenous uncertainty. In each iteration, the algorithm can invoke a simulator to provide a new sample-path. This feature is particularly appealing when a priori representation of uncertainty using scenario trees is either cumbersome or inadequate due to computational and/or timeliness constraints. Such optimization problems are commonly encountered in the operations of power systems with significant renewable penetration. Due to the intermittent nature of renewable resources such as wind and solar, a scenario tree representation may be difficult (perhaps even impossible) to create within the timeliness constraints. State-of-the-art numerical weather prediction and other time series models are known to be more accurate descriptors of such uncertainty. Therefore, optimization algorithms which use sample-paths simulated from such models yield more reliable plans and cost estimates \cite{Gangammanavar2018, Gangammanavar2016}.}

\subsection{Subgradient and Incumbent Selection} \label{sect:policies}
In this section we address two important components of the SDLP algorithm: the ``argmax'' procedure to identify the subgradient of a SDLP approximation at non-root stage that is used during the backward recursion, and the selection of an incumbent solution for the proximal term used during timestaged decision simulation. 

\subsubsection{Subgradient Selection}\label{sect:subgradientPolicy} During the backward recursion, we build a lower bound to the sample average function $H_{t-}^k$ using the best lower bounding affine functions from the collection $\minorants{t}{k}$ for all $\omega_{t} \in \Omega_t^k$. This procedure is accomplished differently based on whether the observation belongs to the current sample-path $\future{\omega}{0}^k$, or not. We utilize the collection of dual vertices $\Pi_t^k$ identified during the course of the algorithm for this purpose. We denote by $i(\pi_t)$ the iteration in which the dual vertex $\pi_t \in \Pi_t^k$ was generated. As seen in \eqref{eq:sdlpDualSet}, the dual vertex $\pi_{t} \in \Pi_t^k$ depends on $H_t^{i(\pi_t)}$, the sample average function in iteration $i(\pi_t)$. This dependence is reflected in the calculation of coefficients $(\bar{\alpha}_{t+}^{i(\pi_t)}, \bar{\beta}_{t+}^{i(\pi_t)})$ and the term $\bar{\rho}_{t+}^{i(\pi_t)}$ that defines the feasible set associated with $\pi_t$ (see \eqref{eq:samplMeanDual}).

{\bf For observation $\boldsymbol{\omega_{t}^k}$:} This observation is encountered at stage $t$ along the current sample-path. Consequently in the current backward recursion, we built and solved a SDA$_t^k$ to optimality using $s_{t}^k$ as input. Using the optimal dual solution thus obtained, we compute the coefficients in \eqref{eq:coefft} for the hyperplanes $\ell_{t}^k(s_{t})$ to SDA$_t^k$ at the candidate state. Similar calculations with $\hat{s}_{t}^k$ as input yield the hyperplane $\hat{\ell}_{t}^k(s_t)$ to SDA$_t^k$ at incumbent state $\hat{s}_t^k$.

{\bf For observations $\boldsymbol{\omega_{t} \in \Omega_{t}^k\setminus \{\omega_{t}^k\}}$:} These are the observations not included in the current sample-path, and therefore, no backward recursion optimization is carried out for these observations. Instead, we use an ``argmax'' procedure to identify the subgradient approximations. These subgradients correspond to the best lower bounding affine functions of SDA$_t^k$ for these observations. In order to accomplish this, we maintain a set of dual solutions $\Pi_t^k$ obtained by solving the SDA$_t^i$ in iterations $i \leq k$ as in the case of 2-SD. For each $\omega_t \in \Omega_t^k\setminus \{\omega_t^k\}$, we setup $s_t = (x_t,\omega_t)$, where $x_t$ is computed with $(x_{t-1}^k, \omega_t, u_{t-1}^k)$ as input in \eqref{eq:stDyn}, and identify a dual solution:
\begin{align*}	\pi_t^k(\omega_t) \in \argmax \bigg \{\bigg(\frac{i(\pi_t)}{k}\bigg)^{T-t} \inner{\pi_t}{(b_t - C_tx_t)}~|~\pi_t \in \Pi_t^k \bigg \}.
\end{align*}
The scaling factor used in the above calculation reflects the scaling of affine functions discussed in Theorem \ref{thm:minorantAsymptotics}. Notice that the set of dual vertices $\Pi_t^k$ changes with iteration which may lead to computational difficulties. We address this issue by using the constancy of the basis index sets that generate these dual vertices. Further discussion of this issue is provided in section \S\ref{sect:incumbentSelection}. Using the dual solution obtained by the above procedure, we can compute the coefficients:
\begin{align*}
	\alpha_t^k(\omega_t) =~& \bigg(\frac{i(\pi_t^k(\omega_t))}{k}\bigg)^{T-t}  [\inner{\pi_t^k(\omega_t)}{b_t} + \bar{\alpha}_{t+}^{~i(\pi_t^k(\omega_t))} ], \\
	\beta_t^k(\omega_t) =~& \bigg(\frac{i(\pi_t^k(\omega_t))}{k}\bigg)^{T-t} [\inner{-C_t}{\pi_t^k(\omega_t)} + \bar{\beta}_{t+}^{~i(\pi_t^k(\omega_t))} ].
\end{align*}
In essence, the above procedure identifies a dual solution $\pi_t^k(\omega_t)$ which was obtained using a SDA$_{t}^{i(\pi_t^k(\omega_t))}$, and scales it appropriately to provide the best lower bounding approximation to the current SDA$_{t}^k$. 

\subsubsection{Incumbent Selection} \label{sect:incumbentSelection}
The procedure described here identifies an incumbent solution at all non-root, non-terminal stages is motivated by the optimal basis propagation policy presented in \cite{Casey2005}. This identification, which is performed during the prediction pass, relies on the basis of  the stage dual approximation (SDA$_t^k$) that appears on the right-hand side of \eqref{eq:samplMeanDual}. To facilitate the discussion here, we have restated SDA$_t^k$ below:
\begin{align} \label{eq:sda_lp}
\max~ \inner{\pi_t}{(b_t - C_tx_t)} \text{ subject to } \inner{D_t}{\pi_t} \leq \bar{\rho}_t^k,~ \pi_t \leq 0, 
\end{align}
where $\bar{\rho}_t^k$ is defined in the expressions following \eqref{eq:samplMeanDual}. In each iteration, the above linear program is solved to optimality along the iteration sample-path and potentially a new basis is discovered. Let $\basis_t^k$ denote the index set whose elements are the rows which are active in \eqref{eq:sda_lp}. Denote by $D_{t,\basis_t^k}$ the submatrix of $D_t$ formed by columns indexed by $\basis_t^k$ (the basis matrix). From standard linear programming results we have that a feasible point is an extreme point of the feasible set if and only if there exists an index set that satisfies
	$\inner{D_{t,\basis_t^k}}{\pi_t^k} = \bar{\rho}_{t,\basis_t^j}^k$.
This index set is added to the collection of previously discovered index sets, that is: $\basisSet_t^k \leftarrow \basisSet_t^{k-1} \cup \basis_t^k$. We use this collection of index sets to construct dual solutions of the linear program in \eqref{eq:sda_lp}. Assumption \assumRef{assum:completeResource} ensures that the optimal set of the dual linear program is non-empty which implies that there exists an index set $\basis_t^j \in \basisSet_t^k$ such that for any arbitrary input state $s_t$ we can write:
\begin{align}\label{eq:incumbGen}
	\hat{u}_{t,i} = D_{t,\basis_t^j}^{-1}(b_t - C_tx_t),~ i \in \basis_t^j; \qquad \hat{u}_{t,i}^j = 0,~i \notin \basis_t^j.
\end{align}
This operation can be written as $\hat{u}_t = R_{\basis_t^j}(b_t - C_tx_t)$, where $R_{\basis_t^j}$ is an $m_t \times n_t$ matrix with rows $[R_{\basis_t^k}]_i = [(D_{t,\basis_t^j})^{-1}]_i$ for $i \in \basis_t^j$ and $[R_{\basis_t^j}]_i = \mathbf{0}$ (a zero vector of length $m_t$) for $i \notin \basis_t^j$. Note that, if $\hat{u}_t^j$ satisfies the constraints of dual of \eqref{eq:sda_lp} then it is a suboptimal basic feasible solution to the dual problem (and if complementarity conditions are also satisfied then it is an optimal solution). We use $\widehat{\mathcal{U}}_t^k(s_t) \subseteq \mathcal{U}_t(s_t)$ to denote the set of basic feasible solutions generated using \eqref{eq:incumbGen} for all index sets in $\mathcal{B}_t^k$. Since \eqref{eq:sda_lp} corresponds to SDA$_t^k$, its dual feasible solutions are feasible to the stage optimization problem \eqref{eq:mslpt}. Using these index sets we define the mapping used for incumbent selection at non-root stages as follows:
\begin{align}\label{eq:incumbMapping}
	\mathcal{M}_t^k(s_t) = \argmin \{ f_t^{k-1}(s_t, \hat{u}_t^j)~|~ \hat{u}_t^j \in \widehat{\mathcal{U}}_t^k(s_t) \} \qquad \forall t \in \mathcal{T}\setminus\{0\}.
\end{align}
We refer to the above mapping as the \emph{basic feasible policy} (BFP) of the MSLP problem. In case the argument that minimizes the right-hand of \eqref{eq:incumbMapping} is not unique, we choose an index set with the smaller iteration index $k$. Notice that the dual LP of \eqref{eq:sda_lp} has cost coefficients that vary over iterations, akin to 2-SD with random cost coefficients in the second-stage \cite{Gangammanavar2020sd}. The steps involved in identifying the BFP, particularly computation of dual solutions in \eqref{eq:incumbGen} and establishing their feasibility, can be implemented in a computationally efficient manner using a sparsity preserving representation of dual solutions. We refer the reader to \cite{Gangammanavar2020sd} for a detailed discussion of this representation and its implementation.

At the root-stage it suffices to maintain a single incumbent solution. This incumbent solution is updated based on predicted objective value reduction at the root-stage:
\begin{align}\label{eq:incumbUpdtt}
	f_0^k(s_0, u_0^k) -  f_0^k(s_0, \hat{u}_0^{k-1})~\leq~q~ [f_0^{k-1}(s_0, u_0^k) - f_0^{k-1}(s_0, \hat{u}_0^{k-1})],
\end{align}
where $q \in (0,1)$ is a given parameter. If the above inequality is satisfied, then the candidate solution at the root node will replace the incumbent solution $\hat{u}_0^{k-1}$ and will serve as the next incumbent solution; that is, $\hat{u}_{0}^k \leftarrow u_{0}^k$ for all $t^\prime \geq t$. On the other hand, if the inequality is not satisfied, then the current incumbent solution for stage $t$ is retained ($\hat{u}_0^k \leftarrow \hat{u}_0^{k-1}$). This update rule is similar to incumbent updates carried out in non-smooth optimization methods including regularized 2-SD \cite{Higle1994, Higle1999}.  

\section{Convergence Analysis}\label{sect:convergenceAnalysis}
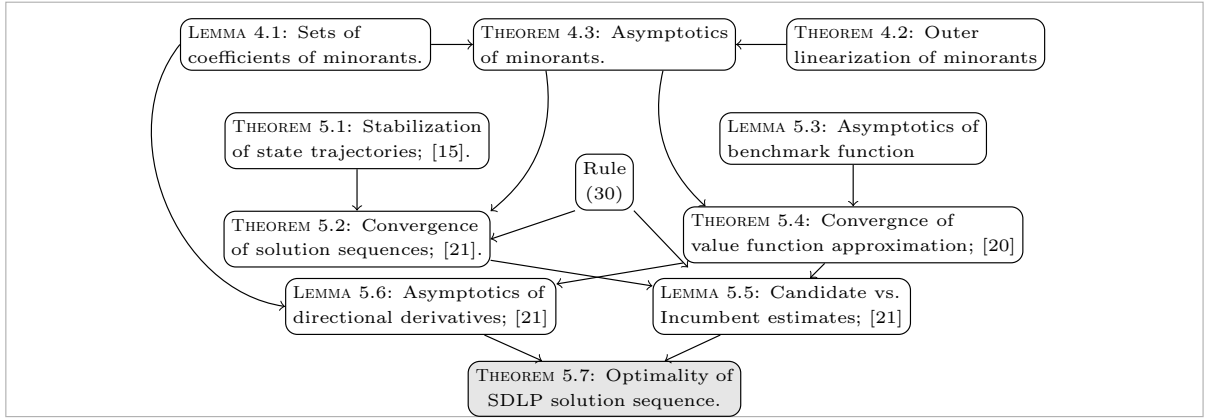
\begin{figure}
	\centering
	\fontsize{6}{7.2}\selectfont
\begin{tikzpicture}[->,blk/.style={rectangle,rounded corners,draw,align=left}]
	\node[draw=black!30, minimum height=5cm, minimum width=0.9\textwidth] (bg) at (0,0) {};

	\node[blk,below=0.2cm of bg.north] (mino) {\theoremRef{thm:minorantAsymptotics}: Asymptotics \\of minorants.};
	\node[blk,left=0.5cm of mino] (coeff) {\lemmaRef{lemma:coeff}: Sets of \\coefficients of minorants.};	
	\node[blk,right=0.6cm of mino] (lb) {\theoremRef{thm:outerLinearization}: Outer \\linearization of minorants};

	\node[blk,below left=0.5cm and -0.2cm of mino] (accumPts) {\theoremRef{thm:converge_predictu}: Stabilization\\ of state trajectories; \cite{Gangammanavar2020sd}.};	
	\node[blk,below=0.5cm of accumPts] (incumbAccum_t) {\theoremRef{thm:converge_optu}: Convergence \\ of solution sequences; \cite{Higle1994}.};	
	
	\node[blk,below right=0.5cm and -0.2cm of mino] (fnConverge) {\lemmaRef{lemma:fnConverge}: Asymptotics of \\benchmark function};	
	\node[blk,below=0.5cm of fnConverge] (valFnConvergence) {\theoremRef{thm:valFnConvergence}: Convergnce of \\value function approximation; \cite{Higle1991}};	
	
	\node[blk,below left=2.5cm and -1cm of mino] (dStability) {\lemmaRef{lemma:directionalStability}: Asymptotics of\\ directional derivatives; \cite{Higle1994}};
	\node[blk,below right=2.5cm and -1cm of mino] (c_vs_i) {\lemmaRef{lemma:candid_v_incumbEst}: Candidate vs.\\ Incumbent estimates; \cite{Higle1994}};	

	\node[blk,below=3.5cm of mino, fill=gray!20,align=center] (opt) {\theoremRef{thm:optimality}: Optimality of\\ SDLP solution sequence.};

	\node[blk,below=1cm of mino,align=center] (incumb) {Rule \\ \eqref{eq:incumbUpdtt}};	
	
	\draw (coeff) to (mino); 
	\draw (lb) to (mino);	
	\draw[bend right=60] (coeff.west) to (dStability.west); 	
	\draw[bend left] (mino.south)++(-0.7,0) to (incumbAccum_t); 
	\draw[bend right] (mino.south)++(0.7,0) to (valFnConvergence);
	
	\draw (fnConverge) to (valFnConvergence); 
	\draw (accumPts) to (incumbAccum_t); 
	\draw (incumb.south east) to (1,-0.7)++(c_vs_i);	
	\draw (incumb.south west) to (incumbAccum_t.east);	

	\draw (valFnConvergence) to (c_vs_i);
	\draw (incumbAccum_t) to (c_vs_i);
	\draw (valFnConvergence) to (dStability);
		
	\draw (c_vs_i) to (opt);	
	\draw (dStability) to (opt);	
	
\end{tikzpicture}
	\caption{Sketch of SDLP Analysis}\label{fig:sdlpAnalysis}%
\end{figure}

In this section, we present the convergence results for SDLP. We begin by discussing the behavior of the sequence of states and decisions generated by the SDLP algorithm, then proceed to show the convergence of value function estimates. Finally, we show that the incumbent solution sequence at the root-stage $\{\hat{u}_0^k\}$ converges and establish the optimality of the accumulation point. The SDLP convergence analysis is built upon the results of the 2-SD algorithm \cite{Higle1991}, its regularized variant \cite{Higle1994}, and 2-SD for 2-SLPs with random cost coefficients \cite{Gangammanavar2020sd}. The Fig. \ref{fig:sdlpAnalysis} illustrates the development of the SDLP convergence analysis. The cited references serve as pointers to related results in the two-stage setting. 

\subsection*{State and decision accumulation points}
Under assumption \assumRef{assum:indep}, we have a finite number possible sample-paths over the horizon. We use $\mathcal{P}_t$ to denote the set of all sample-paths until stage $t$. We focus on the evolution of states and decisions along these sample-paths. 

\begin{theorem}\label{thm:converge_predictu}
Suppose assumptions \assumRef{assum:compact}-\assumRef{assum:indep} hold. Let $\{\hat{u}_0^k\} \subseteq \mathcal{U}_0$ denote any infinite sequence of root-stage incumbent solutions. There exists a subsequence $\mathcal{K}_0$ of iterations such that $\{\hat{u}_0^k\}_{\mathcal{K}_0}$ has an accumulation point. In subsequent stages, for all possible paths $\history{\omega}{t} \in \mathcal{P}_t$ there exists a subsequence of iterations indexed by $\mathcal{K}_t(\history{\omega}{t})$ such that the sequence $\{\hat{u}_{t}^k(\hat{s}_t^k)\}_{k \in \mathcal{K}_t(\history{\omega}{t})}$ has an accumulation point.
\end{theorem}
\begin{proof}
Consider the optimization problem on the right-hand side of  $\eqref{eq:samplMeanDual}$ for a given $t$ in its dual form:
\begin{align}
	\min~\{\inner{\bar{\rho}_t^k}{u_t} ~|~ D_t u_t \leq b_t - C_tx_t, u_t \geq 0\}. \notag
\end{align}
Recall that the feasible set of the above problem is denoted as $\mathcal{U}_t(s_t)$. Let $\mathbb{D}(u_t,s_t) := \argmin \{||u_t - u||^2, u \in \mathcal{U}(s_t) \}$. A slight variant of Hoffman's lemma (see Lemma \ref{lemma:hoffmanVariant} in the appendix) leads us to conclude that for any $s_t, s_t^* \in \text{dom}~ \mathcal{U}_t$ and any $u_t \in \mathcal{U}_t(s_t)$ that 
	$\mathbb{D}(u_t,s_t) \leq \gamma \|(b_t - C_tx_t) - (b_t - C_tx_t^*)\|$.
Here, $\gamma > 0$ is the Lipschitz constant of the mapping $\mathbb{D}(\cdot)$ which depends only on the recourse matrix $D_t$. In other words, the feasible set $\mathcal{U}_t(\cdot)$ is Lipschitz continuous in the above sense. It follows that it is possible to choose an extreme point $\hat{u}_t(s_t) \in \mathcal{U}_t(s_t)$ such that $\hat{u}_t(s_t)$ is continuous on $\text{dom}~\mathcal{U}_t$. Moreover, the polyhedral set $\mathcal{U}_t$ has a finite number of extreme points. Therefore, the BFP outlined in \S \ref{sect:incumbentSelection} is a continuous piecewise linear mapping. 

For the root-node the feasible set $\mathcal{U}_0$ is compact by \assumRef{assum:compact}, hence there exists a subsequence of iterations indexed by $\overline{\mathcal{K}}_0$ such that $\{\hat{u}_0^k\}_{k \in \overline{\mathcal{K}}_0} \rightarrow \bar{u}_0$. Following \assumRef{assum:indep}, there exists an infinite subsequence $\mathcal{K}_1(\history{s}{1}) \subseteq \overline{\mathcal{K}}_0$ such that the algorithm selects sample-path $\history{\omega}{1} \in \mathcal{P}_1$. Since $\{\hat{u}_0^k\}_{k \in \overline{\mathcal{K}}_0}$ converges and $x_0$ is fixed, the sequence of endogenous state $\{x_1^k\}_{k \in \mathcal{K}_1(\history{s}{1})}$ converges to $\bar{x}_1(\history{s}{1})$. For the sample-path $\history{\omega}{1}$, since the sequence of input states $\{x_1^k\}_{k \in \mathcal{K}_1(\history{s}{1})}$ converges, the continuity of BFP implies that the corresponding sequence of incumbent solutions $\{\hat{u}_1^k(\hat{s}_1^k)\}_{k \in \mathcal{K}_1(\history{s}{1})}$ has a converging subsequence. Let $\overline{\mathcal{K}}_1(\history{s}{1})$ denote this subsequence. Therefore, we have $\{\hat{u}_1^k(\hat{s}_1^k)\}_{k \in \overline{\mathcal{K}}_1(\history{s}{1})} \rightarrow \bar{u}_t(\history{s}{1})$.

Now consider an arbitrary stage $t > 1$. For any sample-path $\history{\omega}{t} \in \mathcal{P}_{t}$, once again assumption \assumRef{assum:indep} guarantees that there exists an infinite subsequence of $\mathcal{K}_{t}(\history{s}{t}) \subseteq \mathcal{K}_{t-}(\history{s}{t-})$ when sample-path $\history{\omega}{t}$ is encountered. Here $\history{\omega}{t} = (\history{\omega}{t-}, \omega_t)$, i.e., sample-path $\history{\omega}{t}$ shares the same observations with $\history{\omega}{t-}$ until stage $t-$. Over this subsequence, the convergence of endogenous state sequence $ \{\hat{x}_{t}^k = \mathcal{D}_{t}(\hat{x}_{t-}^k, \omega_{t}, \hat{u}_{t-}^k)\}_{\mathcal{K}_{t}(\history{s}{t})} \rightarrow \bar{x}_{t}(\history{s}{t})$ ensures the convergence of the incumbent states  $\{\hat{s}_{t}^k\}_{\mathcal{K}_{t}(\history{s}{t})}$. Further, the continuity of BFP applied at stage $t$ ensures that the corresponding sequence of incumbent solutions $\{\hat{u}_t^k(\hat{s}_t^k)\}$ have a converging subsequence. That is, there exists $\overline{\mathcal{K}}_{t}(\history{s}{t}) \subset \mathcal{K}_{t}(\history{s}{t})$ such that $\{\hat{u}_t^k(\hat{s}_t^k)\}_{\overline{\mathcal{K}}_{t}(\history{s}{t})} \rightarrow \bar{u}_t(\history{s}{t})$. Proceeding recursively to the rest of the stages, we conclude the validity of the theorem.
\end{proof}

The above result captures the impact of using the argmin mapping in \eqref{eq:incumbMapping} over a sequence of converging first-stage decisions. A converging sequence results in perturbed stage problems with linear constraints in subsequent stages. A central argument in the above proof relies upon the local Lipschitz continuity of the argmin mapping. Such mappings have previously been studied in \cite{Wets2003lipschitz}. We refer the reader to this reference for a more thorough treatment of inf-projections and the argmin mapping for non-linear optimization problems with linear constraints.

To facilitate the presentation in the remainder of this section, let $\future{\mathcal{P}}{t+}^k \in \child_{t+}^k \times \ldots \times \child_T^k$ denote the set of all possible scenarios from stage-$(t+1)$ to the end of horizon which traverse through observations encountered by the algorithm in the first $k$ iterations. Note that $\future{\mathcal{P}}{t+}^k$ represents the set of possible paths in the future and should not be confused with $\mathcal{P}_t^k$ which represents the set of traversed paths. Stagewise independence allows us to compute the probability estimate of a sample-path $\future{\omega}{t+}^j \in \future{\mathcal{P}}{t+}^k$ as product of frequencies associated with observations along that sample-path, i.e. $p^k(\future{\omega}{t+}^j) = p^k(\omega_{t+1}^j) \times \ldots \times p^k(\omega_{T}^j)$. Let $\future{x}{t+}^j$ and $\future{u}{t+}^j$ denote endogenous state and decision vector, respectively, associated with sample-path $\future{\omega}{t+}^j$. While \theoremRef{thm:converge_predictu} captured the behavior of solutions generated using the incumbent mapping in \eqref{eq:incumbGen} during prediction pass, the next result captures the behavior of solutions generated in optimization pass of the algorithm.

\begin{theorem} \label{thm:converge_optu}
	Suppose assumptions \assumRef{assum:compact} - \assumRef{assum:indep} hold, and $\sigma \geq 1$. Then there exists $\bar{u}_0 \in \mathcal{U}_0(s_0)$ such that the sequence of root-node incumbent decisions generated by the algorithm satisfy $\{\hat{u}_0^k\} \rightarrow \bar{u}_0$. Moreover in every subsequent stage, there exists $\bar{u}_{t}(\history{s}{t}) \in \mathcal{U}_t(\bar{s}_{t}(\history{s}{t}))$ which satisfy dynamics in \eqref{eq:stDyn} and the sequence of solutions generated by the algorithm $\{u_{t}^k(\history{s}{t})\}_{\mathcal{K}_t(\history{s}{t})} \rightarrow \bar{u}_{t}(\history{s}{t})$ for all paths $\history{\omega}{t} \in \mathcal{P}_t$. 
\end{theorem}
\begin{proof}
The proof for the root-stage follows that of regularized master in 2-SD (Theorem 5, \cite{Higle1994}) and the root-node of MSD algorithm (\cite{Sen2014}). Here we present the main parts of the proof and refer the reader to earlier works for detailed exposition. If the incumbent solution $\hat{u}_0^k$ changes infinitely many times, then the optimality condition for regularized approximation (see equation (5) in \cite{Higle1994}) and our choice of $\sigma \geq 1$ suggests that for any candidate solution $u_0^k$ the following holds:
\begin{align}\label{eq:regOpt}
 	f_0^{k-1}(s_0, u_0^k) - f_0^{k-1}(s_0, \hat{u}_0^{k-1}) \leq -\|u_0^k - \hat{u}_0^{k-1}\|^2 \leq 0 \qquad \forall k \geq 1.
\end{align}
In particular, the above condition holds at the iterations when the incumbent was updated by assigning the candidate solution as the new incumbent solution, i.e. $\hat{u}_0^k = u_0^k$. Let $\{k_1,k_2,\ldots,k_m\} \in \mathcal{K}_0$ denote the set of $m$ successive iterations when the incumbent solution was updated starting with an incumbent $\hat{u}_0^{k_0}$. Note that, for any $k_n \in \mathcal{K}_0$, $\hat{u}_0^{k_n-1} = \hat{u}_0^{k_{n-1}}$. Denote by $\Delta^{k_n} := f_0^{k_n-1}(s_0, \hat{u}_0^{k_n}) - f_0^{k_n-1}(s_0, \hat{u}_0^{k_{n-1}})$. Using \eqref{eq:regOpt} over these $m$ updates, we have
\begin{align*}
	\frac{1}{m} \sum_{l=1}^m \Delta^{k_n} = &\frac{1}{m} \sum_{l=1}^m [f_0^{k_n-1}(s_0, \hat{u}_0^{k_n}) - f_0^{k_n-1}(s_0, \hat{u}_0^{k_{n-1}})] \\
	= & \frac{1}{m} \underbrace{[f_0^{k_m-1}(s_0, \hat{u}_0^{k_m}) - f_0^{k_1-1}(s_0, \hat{u}_0^{k_0})]}_{(a)} + \\ &\qquad \frac{1}{m} \sum_{n=1}^m \underbrace{[f_0^{k_n-1}(s_0, \hat{u}_0^{k_n}) - f_0^{k_{n+1}-1}(s_0, \hat{u}_0^{k_n})]}_{(b)}.
\end{align*}
The boundedness of functions $\{f_0^k\}$ implies that (a) above approaches zero, as $m \rightarrow \infty$, and their uniform convergence (\theoremRef{thm:minorantAsymptotics}) implies that (b) converges to zero. Hence, 
\begin{align}\label{eq:diminishingError}
    \lim_{m \rightarrow \infty} \frac{1}{m} \sum_{l=1}^m \Delta^{k_n} = 0,
\end{align}
with probability one. Further, the above result, along with \eqref{eq:regOpt} implies that $\lim_{m \rightarrow \infty} \frac{1}{m}
\sum_{n=1}^m \|\hat{u}_0^{k_n}-\hat{u}_0^{k_{n-1}}\|^2 = 0$. Therefore, we conclude that the sequence of root-node incumbent solutions converges to $\bar{u}_0 \in \mathcal{U}_0$.

At non-root stages, the incumbent solutions are selected using the BFP described in \S\ref{sect:incumbentSelection}. The BFP is built using the bases of \eqref{eq:samplMeanDual} discovered during the course of the algorithm that are identified by the collection of index sets $\mathcal{B}_t$. Since there is a finite collection $\mathcal{B}_t$ of index sets, there exists iteration count $K_t$ large enough such that $\mathcal{B}_t^{k^\prime} = \mathcal{B}_t$ for all $k^\prime \geq K_t$. Let us consider $k > \max_t K_t$ when all the index sets for all non-root, non-terminal stages have been discovered. In these iterations, the procedure in \S\ref{sect:incumbentSelection} results in an incumbent solution such that:
\begin{align*}
	\hat{u}_t^k(\hat{s}_t^k) = \mathcal{M}_t^k(\hat{s}_t^k) \in \argmin_{u_t \in \mathcal{U}_t(\hat{s}_t^k)} ~ \inner{d_t}{u_t} + \sum_{\omega_{t+} \in \Omega_{t+}^{k-1}} p^{k-1}(\omega_{t+}) h_{t+}^{k-1}(\dynamics{t}{\hat{x}_t^k,\omega_{t+},u_t}, \omega_{t+}).
\end{align*}
Consequently, the value associated with $\hat{u}_t^k(\hat{s}_t^k)$ is $H_t^{k-1}(\hat{s}_t^k)$ (see \eqref{eq:samplMean}). The forward pass optimal value associated with the candidate solution $u_t^k(\hat{s}_t^k)$ differs from $H_t^{k-1}(\hat{s}_t^k)$ only the quadratic term. Therefore, we have $H_t^{k-1}(\hat{s}_t^k) \leq F_t^{k-1}(\hat{s}_t^k)$ that can be restated using \eqref{eq:objUpdtt} as:
\begin{align*}
	f_t^{k-1}(\hat{s}_t^k, \hat{u}_t^k(\hat{s}_t^k)) - f_t^{k-1}(\hat{s}_t^k, u_t(\hat{s}_t^k)) \leq 0,
\end{align*}
where $u_t(\hat{s}_t^k)$ is the solution obtained by optimizing the regularized problem used during forward recursion. The quadratic programming optimality conditions of this regularized problem allow us to write the following inequality:
\begin{align*}
	f_t^{k-1}(\hat{s}_t^k, u_t(\hat{s}_t^k)) - f_t^{k-1}(\hat{s}_t^k, \hat{u}_t^k(\hat{s}_t^k)) \leq 0.
\end{align*}
The two preceding inequalities together with  (\eqref{eq:regOpt}  for stage $t$) implies that $\|u_t(\hat{s}_t^k) - \hat{u}_t^k(\hat{s}_t^k)\|^2 = 0$. For a sample-path $\history{\omega}{t} \in \mathcal{P}_t$, let $\mathcal{K}_t(\history{s}{t})$ denote the subsequence constructed in the proof of \theoremRef{thm:converge_predictu}. Over this subsequence, the result of \theoremRef{thm:converge_predictu} shows the existence of an accumulation point of $\{\hat{u}_t^k(s_t^k)\}_{k \in \mathcal{K}_t(\history{s}{t})}$, and subsequently, an accumulation point $\bar{u}_t(\history{s}{t})$ of $\{u_t^k(s_t^k)\}_{k \in \mathcal{K}_t(\history{s}{t})}$. Applying the argument to all sample-paths in $\history{\omega}{t} \in \mathcal{P}_t$ completes the proof.
\end{proof}

The limit in \eqref{eq:diminishingError} plays a critical role in showing the existence of an optimal accumulation point of incumbent solutions at the root-stage. Notice that the limit holds when the incumbent is updated infinitely often, i.e., $m \rightarrow \infty$. On the other hand, if the incumbent solution is updated only a finite number of times, then there exists a $K < \infty$ such that $\hat{u}_0^k = \bar{u} \in \mathcal{U}_0$, for all $k > K$. In this case, the optimality of $\bar{u}_0$ is attained only if $\Delta^k \rightarrow 0$. Before we present the optimality of solution sequence, we present the convergence of the value function estimates.

\subsection*{Convergence of Value Function Estimates}
Since our algorithm uses sequential sampling, path-wise forward and backward recursion updates, estimates of probability and sampled minorants we use benchmark functions to verify optimality of the value functions and solutions obtained from them. We next present the construction of these benchmark functions. Note that these function are not computed during the course of the algorithm and are intended only for the purpose of analysis.

For a given input $s_t$, the following is an extensive formulation of the cost-to-function:
\begin{align} \label{eq:pathMean}
	\mathcal{H}_t^k(s_t) =~ & \inner{c_t}{x_t} + \\
	 \min~& \inner{d_t}{u_t} + \sum_{j \in \future{\mathcal{P}}{t+}^k} p^k(\future{\omega}{t+}^j)  \times [\inner{\future{c}{t+}}{\future{x}{t+}^j} + \inner{\future{d}{t+}}{\future{u}{t+}^j}] \notag \\
	& s.t.~u_t \in \mathcal{U}_t(u_0, s_t),~\{u_{t^\prime}^j \in \mathcal{U}_{t^\prime}(u_0,s_{t^\prime}^j)\}_{t^\prime > t} \text{ and non-anticipative}, \notag \\
	& ~~~~~\{x_{t^\prime+}^j = \mathcal{D}_{t^\prime+}(x_{t^\prime}^j,\omega_{t^\prime+}^j,u_{t^\prime}^j)\}_{t^\prime\geq t} .  \notag
\end{align}
In the above formulation, dynamics and non-anticipativity are satisfied starting at stage $t$, and are relative to input $s_t$. This sample average function $\mathcal{H}_t^k$ represents the value associated with input $s_t$ for the remainder of horizon with respect to current observations $\{\child_{i}^k\}_{i=t+}^k$. In order to simplify notation, the dependence of the sample average function on the set $\future{\mathcal{P}}{t+}^k$ is conveyed through the index $k$ in $\mathcal{H}_t^k(s_t)$, as opposed to the more complete $\mathcal{H}_t^k(s_t|\future{\mathcal{P}}{t+}^k)$. 

During forward recursion decisions, $\{u_t\}$ are simulated using approximation $f_t^{k-1}$ in \eqref{eq:forwardt} along the observations dictated by sampling, and during the backward recursion the approximations using subgradients observed along the same sample-path. Next we relate the objective function values encountered during forward and backward recursions. In order to do this, we define $u_t(s_t)$ to be the optimal solution obtained using \eqref{eq:forwardt} during forward recursion with input $s_t$. The forward recursion objective function value $F_t^{k-1}$ associated with this decision is therefore given by:
\begin{align}\label{eq:forwardcost}
	F_t^{k-1}(s_t) := \inner{c_t}{x_t} + &\inner{d_t}{u_t(s_t)} + \sum_{\omega_{t+} \in \child_{t+}^{k-1}} p^{k-1}(\omega_{t+})~ h_{t+}^{k-1}(s_{t+}^k(\omega_{t+})). \notag 
\end{align}
Here $s_{t+}^k(\omega_{t+}) = \mathcal{D}_{t+}(x_t, \omega_{t+}, u_t(s_t))$. In order to study the asymptotic behavior of our algorithm, we investigate how the functions $\mathcal{H}_t^k$, $F_t^k$ and $h_t^k$ relate in value at limiting states. It is worthwhile to note that the sample average approximation in \eqref{eq:samplMeanPrimal}, the extensive formulation in \eqref{eq:pathMean} and the forward recursion objective value in \eqref{eq:forwardcost} are defined only for non-terminal stages as $H_T^k(s_T) = \mathcal{H}_T^k(s_T) = F_T^k(s_T) = h_T(s_T)$ for terminal stage $\forall k$.

\begin{lemma}\label{lemma:fnConverge}
Suppose Assumptions \assumRef{assum:compact}-\assumRef{assum:indep} hold.
	\begin{enumerate}[label=(\roman*)]	
	\item The sequence of functions $\{F^k_t\}_k$ is equicontinuous and uniformly convergent for all $t$. \label{lemma:fnConverge_1}
	\item The sequence of sample average approximation functions $\{\mathcal{H}_t^k\}_k$ converges uniformly to the value function $h_t(\cdot)$ in \eqref{eq:mslpt} for all $t > 0$, with probability one. \label{lemma:fnConverge_2}
	\end{enumerate}
\end{lemma}
\proof
Under Assumption \ref{assum:completeResource}, we have $F_t^k < \infty$ for all $t \in \mathcal{T}\setminus\{T\}$ and $k \geq 1$. (i) Since a regularized problem \eqref{eq:forwardt} with quadratic proximal parameter is used to identify the sequence of solutions in the forward recursion of the algorithm, the optimality conditions of affinely constrained quadratic programs indicate that the solutions $u_t(s_t)$ are piecewise linear. Therefore, the sequence $\{F_t^k\}$ is bounded over a compact space and must have a uniform Lipschitz constant. This leads to the conclusion stated in part \ref{lemma:fnConverge_1} of the lemma. Part \ref{lemma:fnConverge_2} follows from Theorem 7.53 in \cite{Shapiro2014}. \hfill
\endproof

Following the above result, we use $\mathcal{H}_t^k$ as a benchmark for assessing optimality of the SDLP algorithm. We first show the convergence of approximations generated during the course of the algorithm to the true value function in the following theorem. In the two-stage setting, the equivalent result appears as Theorem 3 and Corollary 5 in \cite{Higle1991}.

\begin{theorem} \label{thm:valFnConvergence}
Suppose Assumptions \assumRef{assum:compact}-\assumRef{assum:fixed} hold. At any non-terminal stage $t$, if subsequence $\mathcal{K}_t$ is such that $\{\hat{s}_t^k\}_{k \in \mathcal{K}_t} \rightarrow \bar{s}_{t}$, then 
\begin{align}
	\lim_{k \in \mathcal{K}_t} f_t^k(\hat{s}_t^k, \hat{u}_t^k(\hat{s}_t^k)) = \lim_{k \in \mathcal{K}_t} f_t^{k+1} (\hat{s}_t^k, \hat{u}_t^k(\hat{u}_t^k)) = f_t(\bar{s}_t, \bar{u}_t),
\end{align}
with probability one.
\end{theorem}
\begin{proof}
For terminal stage ($t = T$), continuity of linear programming value function implies that $\lim_{k \in \mathcal{K}_T} h_T(\hat{s}_T^k) = h_T(\bar{s}_T(\history{s}{t}))$. Since $F_T^k$, $H_T^k$ and $\mathcal{H}_T^k$ are all equivalent to $h_T$, the above relation trivially holds. Consequently we have, $\lim_{k \rightarrow \mathcal{K}_T} \hat{\ell}_T^k(\hat{s}_T^k) = h_T(\bar{s}_T)$ and $\lim_{k \rightarrow \mathcal{K}_T} \partial \hat{\ell}_T^k(\hat{s}_T^k) \in \partial h_T(\bar{s}_T)$.

For a non-terminal stage, let $k - \tau$ and $k$ be two successive iterations of subsequence $\mathcal{K}_t$. The forward recursion objective function $F_t^{k-1}(s_t)$ and the backward recursion sample average function $H_t^{k-1}(s_t)$ differ only in the proximal term, and hence $H_t^{k-1}(s_t) \leq F_t^{k-1}(s_t)$ for all $s_t \in \mathcal{S}_t$. In the following, we focus on functions evaluated at $\hat{s}_t^k$, and use $m = \omega_{t+}^k$ and $n$ as an index for set $\Omega_{t+}^k$. The forward recursion objective function value at the current input state can be written as:
\begin{align*}
	F_t^{k-1}(\hat{s}_t^k) =~& \inner{c_t}{\hat{x}_t^k} + \inner{d_t}{u_t(\hat{s}_t^k)} + \sum_{n \in \Omega_{t+}^{k-1}} p^{k-1}(n) h_{t+}^{k-1}(\hat{s}_{nt+}^k).
\end{align*}
The optimality of $u_t(\hat{s}_t^k)$ ensures that the objective function value is associated with $u_t(\hat{s}_t^k)$ is lower than any other feasible solution. If we specifically consider the optimal solution of the dual in \eqref{eq:samplMeanDual}, denoted $\tilde{u}_t(\hat{s}_t^k)$, we have
\begin{align*}
	F_t^{k-1}(\hat{s}_t^k) \leq~& \inner{c_t}{\hat{x}_t^k} + \inner{d_t}{\tilde{u}_t(\hat{s}_t^k)} + \sum_{n \in \Omega_{t+}^{k-1}} p^{k-1}(n) h_{t+}^{k-1}(\tilde{s}_{nt+}^k).
\end{align*}
By adding and subtracting the current approximation of future cost, i.e., $\sum_{n \in \Omega_{t+}^k} p^k(n) h_{t+}^k(\hat{s}_{nt+}) = \sum_{n \in \Omega_{t+}^{k-1}\setminus\{m\}} p^k(n)~ h_{t+}^k(\hat{s}_{nt+}) + p^k(m) h_{t+}^k(\hat{s}_{mt+})$ we obtain
\begin{align*}
	&F_t^{k-1}(\hat{s}_t^k) \leq \inner{c_t}{\hat{x}_t^k} + \inner{d_t}{\tilde{u}_t(\hat{s}_t^k)} + \sum_{n \in \Omega_{t+}^k} p^k(n) h_{t+}^k(\hat{s}_{nt+}^k) + \\ &\sum_{n \in \Omega_{t+}^{k-1}} p^{k-1}(n) h_{t+}^{k-1}(\tilde{s}_{nt+}^k) - \bigg[ \sum_{n \in \Omega_{t+}^{k-1}} \bigg (\frac{k-1}{k} \bigg)p^{k-1}(n)~ h_{t+}^k(\hat{s}_{nt+}^k) + \frac{1}{k} h_{t+}^k(\hat{s}_{mt+}^k) \bigg ].
\end{align*}
From the definition of backward recursion sample average approximation in \eqref{eq:samplMeanDual} and the fact that $h_t(s_t) \geq 0$, we have
\begin{align*}
	F_t^{k-1}(\hat{s}_t^k) \leq~& H_t^k(\hat{s}_t^k) + \sum_{n \in \Omega_{t+}^{k-1}} p^{k-1}(n) \bigg[h_{t+}^{k-1}(\tilde{s}_{nt+}^k) - \bigg(\frac{k-1}{k} \bigg) h_{t+}^k(\hat{s}_{nt+}^k)\bigg].
\end{align*}
From Theorem \ref{thm:minorantAsymptotics}, we have $h_{t+}^k \geq \big(\frac{k-1}{k}\big)^{T-t-1}h_{t+}^{k-1}$. This yields
\begin{align*}
    F_t^{k-1}(\hat{s}_t^k) \leq~& H_t^k(\hat{s}_t^k) + \sum_{n \in \Omega_{t+}^{k-1}} p^{k-1}(n) \bigg[h_{t+}^{k-1}(\tilde{s}_{nt+}^k) - \bigg(\frac{k-1}{k} \bigg)^{T-t} h_{t+}^{k-1}(\hat{s}_{nt+}^k)\bigg].
\end{align*}
Let us focus on the terms within the summation on the right hand side of above inequality, i.e., $\Delta_n^k = h_{t+}^{k-1}(\tilde{s}_{nt+}^k) - \big(\frac{k-1}{k}\big)^{T-t}h_{t+}^k(\hat{s}_{nt+}^k)$. Then
\begin{align*}
	\Delta_n^k =~& h_{t+}^{k-1}(\tilde{s}_{nt+}^k) - h_{t+}^{k-1}(\hat{s}_{nt+}^k) + \bigg(1-\bigg(\frac{k-1}{k}\bigg)^{T-t}\bigg)h_{t+}^{k-1}(\hat{s}_{nt+}^k)\\
	\leq~& |h_{t+}^{k-1}(\tilde{s}_{nt+}^k) - h_{t+}^{k-1}(\hat{s}_{nt+}^k)| + \bigg|\bigg(1-\bigg(\frac{k-1}{k}\bigg)^{T-t}\bigg)h_{t+}^{k-1}(\hat{s}_{nt+}^k)\bigg|.
\end{align*}
The second term in the above equates to zero as $k \rightarrow \infty$. Further, since $\{\hat{s}_t^k\}_{k \in \mathcal{K}_t} \rightarrow \bar{s}_{t}(\history{s}{t})$, for every $\delta > 0$ there exists a $K(\delta) \in \mathcal{K}_t$ such that $\|\hat{s}_t^k - \tilde{s}_t^k\| < \delta$ for all $k > K(\delta)$. Using the uniform equicontinuity of $\{h_t^k\}$ (\theoremRef{thm:minorantAsymptotics}), we have $|h_{t+}^{k-1}(\tilde{s}_{nt+}^k) - h_{t+}^{k-1}(\hat{s}_{nt+}^k)| < \epsilon$. Therefore, we can conclude that $\lim_{k \in \mathcal{K}} F_t^{k-1}(\hat{s}_t^k) - H_t^k(\hat{s}_t^k) \leq \epsilon$, for any $\epsilon > 0$.

To show the inequality in the other direction, we use the fact that $H_t^{k}(\hat{s}_t^k) \leq F_t^{k}(\hat{s}_t^k)$ and the uniform convergence of the sequence $\{F_t^k\}$. This gives us $\lim_{k \in \mathcal{K}_t} H_t^k(\hat{s}_t^k) - F_t^{k-1}(\hat{s}_t^k) \leq \epsilon$. Since inequalities hold in both directions for an arbitrary $\epsilon > 0$, we have
\begin{align} \label{eq:squeeze1}
	\lim_{k \in \mathcal{K}_t} |F_t^{k-1}(\hat{s}_t^k) - H_t^k(\hat{s}_t^k)| = 0 ~~~ (w.p.1).
\end{align}

Now consider the benchmark function $\mathcal{H}_t^k(\hat{s}_t^k)$ that is optimal across all possible sample-paths. Optimality of $\mathcal{H}_t^k(\hat{s}_t^k)$, along with the fact that $h_t^k \leq H_t^k$ (\theoremRef{thm:outerLinearization}), we have
\begin{align*}
	\lim_{k \in \mathcal{K}_t} \mathcal{H}_t^k(\hat{s}_t^k) \leq \lim_{k \in \mathcal{K}_t} h_t^k(\hat{s}_t^k) \leq \lim_{k \in \mathcal{K}_t} H_t^k(\hat{s}_t^k) ~~~ (w.p.1),
\end{align*}

Moreover, the forward recursion objective function value satisfies $\lim_{k \in \mathcal{K}_t} F_t^{k-1}(\hat{s}_t^k) \leq \lim_{k \in \mathcal{K}_t} \mathcal{H}_t^k(\hat{s}_t^k)~ (w.p.1)$. Therefore we have
\begin{align} \label{eq:squeeze2}
	\lim_{k \in \mathcal{K}_t} F_t^{k-1}(\hat{s}_t^k) \leq \lim_{k \in \mathcal{K}_t} \mathcal{H}_t^k(\hat{s}_t^k) \leq \lim_{k \in \mathcal{K}_t} h_t^k(\hat{s}_t^k) \leq \lim_{k \in \mathcal{K}_t} H_t^k(\hat{s}_t^k) ~~~ (w.p.1).
\end{align} Using \eqref{eq:squeeze1} in the above relation and the results in \lemmaRef{lemma:fnConverge}, we conclude that the expresstion \eqref{eq:squeeze2} holds with equality, with probability one.

Since $\{\hat{s}_t^k\}_{k \in \mathcal{K}_t} \rightarrow \bar{s}$, the result in \theoremRef{thm:converge_optu} shows the existence a subsequence $\overline{\mathcal{K}}_t$ such that $\{(\hat{s}_t^k, \hat{u}_t^k)\}_{k \in \overline{\mathcal{K}}_t} \rightarrow (\bar{s}_t,\bar{u}_t)$. Using, the uniform convergence of the sequence of minorants $\{h_t^k\}$ and benchmark function $\{\mathcal{H}_t^k\}$ (\theoremRef{thm:minorantAsymptotics} and \lemmaRef{lemma:fnConverge}, respectively), we conclude that the function values  $\{f_t^k(\hat{s}_t^k,\hat{u}_t^k(\hat{s}_t^k))\}$ converge to the optimal value at the accumulating state $\bar{s}_t$, with probability one.
\end{proof}
The above result holds at all non-terminal stages over any converging subsequence of states $\{\hat{s}_t^k\}_{k \in \mathcal{K}_t}$. Under the assumption of finite support \assumRef{assum:indep}, the algorithm will generate such subsequences as illustarted in \theoremRef{thm:converge_predictu} and \theoremRef{thm:converge_optu}. 
\subsection*{Optimality of the Incumbent Solution Sequence}
Before establishing the optimality of the root-stage incumbent solution sequence, we establish the limiting relationship between the value function estimate at the candidate solution $f_0^{k-1}(s_0, u_0^k)$ and estimate at the incumbent solutions $f_0^{k-1}(s_0, \hat{u}_0^{k-1})$.  As a consequence of \ref{thm:valFnConvergence}, the root-stage value function is equivalent to the value function of a 2-SLP. This equivalent 2-SLP has the first-stage cost equal to $\inner{c_0}{x_0} + \inner{d_0}{u_0}$ and the expected recourse value given by $\sum_{\omega_1 \in \Omega_1} p(\omega_1) h_1(s_1(\omega))$. The function $h_1(\cdot)$ is the optimal cost-to-go value starting from stage $1$ which is attained for the limiting states $\bar{s}_1(\omega_1)$ for all $\omega_1 \in \Omega_1$. With this perspective, the following lemma parallels a result from \cite{Higle1994} (Theorem 3). We present the proof for the case when the incumbent changes infinitely often and refer the reader to \cite{Higle1994} for the case when the incumbent changes finitely often.

\begin{lemma}\label{lemma:candid_v_incumbEst}
Let $\{u_0^k\}_{k=1}^\infty$ and $\{\hat{u}_0^k\}_{k=1}^\infty$ denote the sequence of candidate and incumbent solutions identified by SDLP, respectively. With probability one,
\begin{align} \label{eq:candid_v_incumbEst}
	\limsup_{k \rightarrow \infty} f_0^{k-1}(s_0, u_0^k) - f_0^{k-1}(s_0, \hat{u}_0^{k-1}) = 0.
\end{align}
\end{lemma}
\begin{proof}
Let $\{k_n\}_{n \in \mathcal{K}_0}$ represent the sequence of iterations at which the incumbent is changed. If $\mathcal{K}_0$ is an infinite set, then as a consequence of the incumbent update rule \eqref{eq:incumbUpdtt} and \theoremRef{thm:valFnConvergence}, we have
\begin{align*}
	\lim_{m \rightarrow \infty} \frac{1}{m} \sum_{n=1}^m \Delta^{k_n} \leq \limsup_{n \rightarrow \infty} \Delta^{k_n} \leq 0.
\end{align*}
From \eqref{eq:diminishingError}, there exists a subsequence $\mathcal{K}_0^* \subset \mathcal{K}_0$ such that
\begin{align*}
	\lim_{k \in \mathcal{K}_0^*} \Delta^k = 0.
\end{align*}
Since, $\Delta^k = f_0^{k-1}(s_0, \hat{u}_0^k) - f_0^{k-1}(s_0, \hat{u}_0^{k-1})$ and $\hat{u}_0^k = u_0^k$, for all $k \in \mathcal{K}_0^*$, establishes \eqref{eq:candid_v_incumbEst}.
\end{proof}

The following result captures the asymptotic behavior of the directional derivatives of the sequence of first-stage objective function approximations. Specifically, it relates the directional derivatives of the approximate value function to that of the true value function.
\begin{lemma}\label{lemma:directionalStability}
Let $u_t \in \mathcal{U}_t(s_t)$. Define $\delta_t^k(u_t) = \frac{u_t - u_t^{k-1}}{\|u_t - u_t^{k-1}\|}$ and $\bar{\delta}_t(u_t) = \frac{u_t - \bar{u}_t}{\|u_t - \bar{u}_t\|}$. For any sequence $\mathcal{K}$ such that $\{u_t^k\}_{k \in \mathcal{K}} \rightarrow \bar{u}_t$, $\bar{u}_t \in \mathcal{U}_t(s_t)$, then 
\begin{align}
	\lim_{k \in \mathcal{K}} (f_t^k)^\prime(s_t, u_t^{k-1}; \delta_t^k(u_t)) \leq f_t^\prime(s_t, \bar{u}_t; \bar{\delta}_t(u_t)), 
\end{align}
with probability one.
\end{lemma}
\begin{proof}
Since $\{u_t^k\}_{k \in \mathcal{K}} \rightarrow \bar{u}_t$, we have $\delta_t^k(u_t) \rightarrow \bar{\delta}_t(u_t)$, for all $u_t \in \mathcal{U}_t(s_t)$. Note that, 
\begin{align*}
	\beta_{t+}^k(\omega_{t+}) = (h_{t+}^k)^\prime(x_{t+}(\omega_{t+}), \omega_{t+}).
\end{align*}
Following \lemmaRef{lemma:coeff} (b) and \theoremRef{thm:valFnConvergence}, we have $\limsup_{k \in \mathcal{K}} (f_t^k)^\prime(s_t, u_t^{k+1}; \delta_t^k(u_t)$ is finite and  further, there is exists a subsequence $\overline{\mathcal{K}} \subset \mathcal{K}$ such that 
\begin{align} \label{eq:subgradientAccum}
	\{\beta_{t+}^k(\omega_{t+})\}_{k \in \overline{\mathcal{K}}} \rightarrow \bar{\beta}_{t+}(\omega_{t+}) \in \partial h_{t+}(x_{t+}(\omega_{t+}), \omega_{t+}).
\end{align}
Using the definition of $f_t^k$ in \eqref{eq:objUpdtt}, we have
\begin{align*}
	&(f_t^k)^\prime(s_t, u_t^{k+1}; \delta_t^k(u_t)) = \inner{d_t + \sum_{\omega_{t+} \in \child_{t+}^k} p^k(\omega_{t+}) (h_{t+}^k)^\prime(x_{t+}(\omega_{t+}), \omega_{t+})}{\delta_t^k(u_t}.
\end{align*}
This implies that
\begin{align*}
	\lim_{k \in \overline{\mathcal{K}}} (f_t^k)^\prime(s_t, u_t^{k+1}; \delta_t^k(u_t)) =~& \limsup_{k \in \mathcal{K}} (f_t^k)^\prime(s_t, u_t^{k+1}; \delta_t^k(u_t)) \\
	=~& \inner{d_t + \sum_{\omega_{t+} \in \child_{t+}^k} p^k(\omega_{t+}) (h_{t+}^k)^\prime(x_{t+}(\omega_{t+}), \omega_{t+})}{\delta_t^k(u_t)}. 	
\end{align*}
Let $\bar{d}_t = d_t + \expect{\bar{\beta}(\tilde{\omega})}{}$. Using \eqref{eq:subgradientAccum} and the fact that $p^k(\omega) \rightarrow p(\omega)$, almost surely, we have
\begin{align*}
	\limsup_{k \in \mathcal{K}} (f_t^k)^\prime(s_t, u_t^{k+1}; \delta_t^k(u_t)) = ~& \inner{\bar{d}_t}{\bar{\delta}_t(u_t)} \\
	\leq ~& \max\{\inner{v}{\bar{\delta}_t(u_t)}~|~ v \in \partial f_t(s_t,\bar{u}_t; \bar{\delta}(\bar{u}_t)\} \\
	=~& f^\prime(s_t, \bar{u}_t; \bar{\delta}(\bar{u}_t).
\end{align*}
\end{proof}

The above lemma mirrors a similar result  for regularized 2-SD that appeared in \cite{Higle1994} (as Lemma 4). We are now in a position to establish the optimality of the accumulation point of the sequence of root-stage incumbent solutions.
\begin{theorem} \label{thm:optimality}
	Suppose Assumptions \assumRef{assum:compact}-\assumRef{assum:indep} hold and $\underline{\sigma} \geq 1$, then the SDLP algorithm produces a sequence incumbent solutions at root-stage $\{\hat{u}_0^k\} \rightarrow u_0^*$ and $u_0^*$ is optimum, with probability one.
\end{theorem}
\begin{proof}
Using the optimality condition of regularized root-stage problem \eqref{eq:regOpt}, the result in \lemmaRef{lemma:candid_v_incumbEst} implies that there exists a subsequence $\mathcal{K}_0^*$ such that
\begin{align*}
	\lim_{k \in \mathcal{K}_0^*} f_0^{k-1}(s_0, u_0^k) - f_0^{k-1}(s_0, \hat{u}_0^{k-1}) + \|u_0^k - \hat{u}_0^{k-1}\| = 0,
\end{align*}
with probability one. Let $\overline{\mathcal{K}}_0^* \subset \mathcal{K}_0^*$ be such that $\{\hat{u}_0^k\}_{k \in \overline{\mathcal{K}}_0^*} \rightarrow \bar{u}_0$. Let $u_0 \in \mathcal{U}_0$ be such that $u_0 \neq \bar{u}_0$. We define
\begin{align*}
	\delta_0^k(u_0) = \frac{u_0 - u_0^k}{\|u_0 - u_0^k\|},~ \text{and} ~\bar{\delta}_0(u_0) = \frac{u_0 - \bar{u}_0}{\|u_0 - \bar{u}_0\|}.
\end{align*}
Optimality of $u_0^k$ implies that 
\begin{align*}
	& f_0^k(s_0, u_0^k) + \frac{\sigma}{2}\|u_0^k - \hat{u}_0^{k-1}\|^2 \leq f_0^k(s_0, u_0^k + \delta_0^k(u_0)) + \frac{\sigma}{2}\|(u_0^k + \delta_0^k(u_0)) - \hat{u}_0^{k-1}\|^2. \\
	\Rightarrow & [f_0^k(s_0, u_0^k + \delta_0^k(u_0)) - f_0^k(s_0, u_0^k)] + \\
	& \hspace*{4cm}\frac{\sigma}{2}[\|(u_0^k + \delta_0^k(u_0)) - \hat{u}_0^{k-1}\|^2 - \|u_0^k - \hat{u}_0^{k-1}\|^2] \geq 0. \\
	\Rightarrow & (f_0^k)^\prime (s_0, u_0^k; \delta_0^k(u_0)) + \frac{\sigma}{2}[\|(u_0^k + \delta_0^k(u_0)) - \hat{u}_0^{k-1}\|^2 - \|u_0^k - \hat{u}_0^{k-1}\|^2] \geq 0.
\end{align*}
Taking limits along $\overline{\mathcal{K}}_0^*$, the second term equates to zero. Therefore, we have
\begin{align*}
	0 \leq \liminf_{k \in \mathcal{K}_0^*} (f_0^k)^\prime (s_0, u_0^k; \delta_0^k(u_0)) \leq \limsup_{k \in \mathcal{K}_0^*} (f_0^k)^\prime (s_0, u_0^k; \delta_0^k(u_0)) \leq f^\prime(s_0, \bar{u}_0; \bar{\delta}_0(u_0)).
\end{align*} 
The last inequality follows from \lemmaRef{lemma:directionalStability}. Since $f_0(\cdot)$ is convex function and the above statement implies that the directional derivatives $\bar{\delta}_0(u_0)$ are non-negative for an arbitrary $u_0 \in \mathcal{U}_0$. We must have that $\bar{u}_0$ must be an optimal solution.
\end{proof}

\section{Conclusions}\label{sect:conclusions}
The SDLP algorithm extends the regularized 2-SD algorithm \cite{Higle1994} to the MSLP setting where the underlying stochastic process exhibits stagewise independence. The algorithm addresses the state variable formulation of MSLP problems by employing sequential sampling. In this sense, it is a counterpart to the MSD algorithm of \cite{Sen2014} which was designed for a case where the underlying uncertainty has a scenario tree structure. The algorithm presented in this paper incorporates several additional advantages granted by the stagewise independence property. We conclude here by noting the salient features of the SDLP algorithm:
\begin{enumerate}
	\item The algorithm uses a single sample-path both for simulating decisions during the forward recursion and updating approximations during backward recursion. In any iteration, compared to SDDP which requires solving subproblems corresponding to all outcomes at all stages and for all sample-paths simulated during the forward pass, SDLP uses two subproblem solves at each stage. This significantly reduces the computational burden of solving MSLP problems.
	\item The method uses quadratic regularization terms at all non-terminal stages which alleviates the need to retain all the minorants generated. This allows us to retain a finite-sized approximation in all stages, further improving its computational advantage. 
	\item The BFP described in \S\ref{sect:incumbentSelection} is the first to provide a data-driven policy for MSLP. This mapping overcomes the need to store incumbent solutions that, either explicitly or implicitly, depending on the entire history of state evolution, and can be used with other regularized MSLP algorithms. Our convergence results show that the optimality of the accumulation points of a subsequence of incumbent solutions is preserved even when such a mapping is employed.
	\item SDLP incorporates sampling within the optimization step, and thereby, optimizes an SAA with increasing sample size. This feature enables SDLP to solve the MSLP problems to greater accuracy by incorporating additional observations at any stage without having to re-discover the structural information of an instance to build/update the approximations. The adaptive nature allows the algorithm to be terminated upon attaining a desired level of accuracy. This opens the avenue to design statistical optimality rules for multistage setting akin to those developed for 2-SLP \cite{Higle1999, Sen2016}. 
\end{enumerate} 
The computational advantages of SDLP were revealed in our companion paper \cite{Gangammanavar2018}. In that paper, we applied the SDLP algorithm to a MSLP model for distributed storage control in the presence of renewable generation uncertainty. The computational results compare our algorithm with SDDP applied to a SAA of the original model. The sample-paths used to set up the SAA and those used within the SDLP algorithm were simulated using an autoregressive moving-average time series model. The computational results of that paper indicate that SDLP provides solutions that are not only reliable but are also statistically indistinguishable from SDDP, while significantly improving the computational times. The computational advantage of SDLP over SDDP can be attributed to the algorithm design. Namely, (i) the forward and backward recursion calculations are carried along only one sample-path in each iteration, and (ii) the use of regularization helps us maintain a finite sized optimization problem at every non-terminal stage. Note that we are only referring to calculations within any particular iteration. In this sense our comparison is incomplete. However, carrying out a full theoretical comparison of SDLP and SDDP is beyond the scope of this paper. Nevertheless, we point the reader to recent results related to iteration complexity of the SDDP algorithm \cite{Lan2020complexity} and the sublinear rate of convergence for 2-SD in \cite{Liu2020asymptotic}. We plan to undertake the convergence rate analysis, (sample and iteration complexity) of SDLP in our future research endeavors. In any case, the results in \cite{Gangammanavar2018} provide the first evidence of computational benefits provided by a sequential sampling approach in a multistage setting.

\appendix
\section{A Variant of Hoffman's Lemma} \label{append:proofs}
In this appendix we present a variant of the Hoffman's Lemma that is integral to the proof of \theoremRef{thm:converge_predictu}.
\begin{lemma}\label{lemma:hoffmanVariant}
Let $\mathfrak{U}(x,\rho)$ be the set of optimal primal solutions of problem \eqref{eq:sda_lp}. Then there exists a positive constant $\chi$, depending only on $C$ and $D$, such that for any $(x,\rho),(x^\prime,\rho^\prime) \in \text{dom}~\mathfrak{U}$ and any $u \in \mathfrak{U}(x,\rho)$,
\begin{align} \label{eq:optSolContinuousMapping}
	\text{dist}(u,\mathfrak{U}(x^\prime,\rho^\prime)) \leq \chi \|x-x^\prime\|.
\end{align}
\end{lemma}
\begin{proof}
The linear program can be written in an equivalent form: 
\begin{align} \label{eq:stageProblemEqui}
	\min_{\eta \in \mathbb{R}}~ \eta~ \text{subject to}~ Du \leq b - Cx,~\inner{\rho}{u} - \eta \leq 0.
\end{align}
Denote by $\mathcal{E}(x,\rho) := \{(u,\eta)~|~ Du \leq b - Cx,~\inner{\rho}{u} - \eta \leq 0\}$ the set of feasible points of \eqref{eq:stageProblemEqui}. Let $(x,\rho),(x^\prime,\rho^\prime) \in \text{dom}~ \mathfrak{U}$ and consider a point $(u,\eta) \in \mathcal{E}(x,\rho)$. Note that for any $a \in \mathbb{R}^{n}$ we have $\|a\| = \sup_{\|z\|_* \leq 1} \inner{z}{a}$, where $\|\cdot\|_*$ is the dual of the norm $\|\cdot\|$. Using this we have
\begin{subequations}\begin{align}
	\text{dist}((u,\eta), \mathcal{E}(x^\prime,\rho^\prime)) = &\inf_{(u^\prime,\eta^\prime) \in \mathcal{E}(x^\prime,\rho^\prime)} \|(u,\eta) - (u^\prime,\eta^\prime)\| \\
	= & \inf_{\substack{Du^\prime \leq b - Cx^\prime \\ \inner{\rho^\prime}{u^\prime} - \eta^\prime \leq 0}} \sup_{\|(z_0,z_1)\|_* \leq 1} \inner{z_0}{(u - u^\prime)} + z_1(\eta - \eta^\prime) \\
	= & \sup_{\|(z_0,z_1)\|_* \leq 1} \inf_{\substack{Du^\prime \leq b - Cx^\prime \\ \inner{\rho^\prime}{u^\prime} - \eta^\prime \leq 0}} \inner{z_0}{(u - u^\prime)} + z_1(\eta - \eta^\prime).
\end{align}\end{subequations}
The interchange between the minimum and maximum operators follows from Theorem 7.11 in \cite{Shapiro2014}. By using the change of variables $\Delta u = (u - u^\prime)$ and $\Delta \eta = (\eta - \eta^\prime)$, we can rewrite the inner term as:
\begin{align}
	\inf_{\substack{D\Delta u \geq Du - (b - Cx^\prime) \\ \inner{\rho^\prime}{\Delta u} - \Delta \eta \geq \inner{\rho^\prime}{u} - \eta}} \inner{z_0}{\Delta u} + z_1 \Delta \eta.
\end{align}
The dual of the above problem is given by:
\begin{align}
	\sup_{\substack{\lambda \geq 0, \mu \geq 0\\ \inner{D}{\lambda} + \mu\rho^\prime = z_0 \\ \mu = -z_1}} \inner{\lambda}{[Du - (b-Cx^\prime)]} + \mu [\inner{\rho^\prime}{u} - \eta].
\end{align}
It follows from $\|(z_0,z_1)\|_* \leq 1$ that 
\begin{align}
	\text{dist}((u,\eta), \mathcal{E}(x^\prime,\rho^\prime)) = \sup_{\substack{\lambda \geq 0,~ 0 \leq \mu \leq 1\\ \|\inner{D}{\lambda} + \mu\rho^\prime\|_* \leq 1}} \inner{\lambda}{[Du - (b-Cx^\prime)]} + \mu [\inner{\rho^\prime}{u} - \eta].
\end{align}

Using similar arguments used in Lemma \ref{lemma:coeff} \ref{lemma:coeff_b}. we can establish that the cost coefficients of \eqref{eq:sda_lp} generated within SD-based methods for non-terminal stages belong to compact sets. Therefore, we can assume without loss of generality that $\|\rho^\prime\|_* \leq M$. We obtain a relaxation of the above dual problem by replacing the constraint $\|\inner{D}{\lambda} + \mu\rho^\prime\|_* \leq 1$ with the constraint $\|\inner{D}{ \lambda}\|_* \leq 1 + M$. Let $(\hat{\lambda},\hat{\mu})$ be an optimal solution of the relaxed dual problem. We can assume without loss of generality that $\|\cdot\|$ is the $\ell_1$ norm, and hence its dual is the $\ell_\infty$ norm. For such a choice of a polyhedral norm, we have that the feasible set of the relaxed dual problem is polyhedral. Therefore, $\hat{\lambda}$ is the an extreme point of the set $\{\lambda~|~ \|\inner{D}{\lambda} \|_* \leq 1 + M\}$. This implies that $\|\lambda\|_*$ can be bounded by a constant $\chi_1$ which depends only on $D$.

Since $(u,\eta) \in \mathcal{E}(x,\rho)$, and hence $Du - (b-Cx)\leq 0$ and $\inner{\rho}{u}-\eta \leq 0$, we have
\begin{align*}
	\inner{\hat{\lambda}}{[Du - (b-Cx^\prime)]} ~=~ & \inner{\hat{\lambda}}{[Du - (b-Cx)]} + \inner{\hat{\lambda}}{C(x-x^\prime)} \\
									  \leq~ & \inner{\hat{\lambda}}{C(x-x^\prime)} \leq~ \|\lambda\|_* \|C\|\|x-x^\prime\|.
\end{align*}
Further, notice that $\hat{\mu}[\inner{\rho}{u} - \eta] \leq 0$ and $\|C\| \leq \chi_2$. This leads us to conclude that 
\begin{align}
	\text{dist}((u,\eta), \mathcal{E}(x^\prime,\rho^\prime)) = \|\hat{\lambda}\|_*\|C\|\|x-x^\prime\| \leq \chi_1\chi_2\|x-x^\prime\|.
\end{align} This implies that  \eqref{eq:optSolContinuousMapping} is true.
\end{proof}

\bibliographystyle{plain}
\bibliography{sdlpBiblio}

\end{document}